\ifdefined\pdfoutput\pdfoutput=1\fi
\documentclass[11pt]{amsart}
\usepackage[margin=1.15in]{geometry}
\usepackage{amsmath,amssymb,amsthm,mathtools}
\usepackage{booktabs,array}
\usepackage[hidelinks]{hyperref}
\usepackage{tikz}
\usetikzlibrary{calc,arrows.meta,positioning}
\usepackage[expansion=false]{microtype}

\newtheorem{maintheorem}{Theorem}

\newtheorem{maincorollary}[maintheorem]{Corollary}

\newtheorem{theorem}{Theorem}[section]
\newtheorem{proposition}[theorem]{Proposition}
\newtheorem{lemma}[theorem]{Lemma}
\newtheorem{corollary}[theorem]{Corollary}
\theoremstyle{definition}
\newtheorem{definition}[theorem]{Definition}
\newtheorem{question}[theorem]{Question}
\theoremstyle{remark}
\newtheorem{remark}[theorem]{Remark}
\numberwithin{equation}{section}

\newcommand{\Z}{\mathbb{Z}}
\newcommand{\R}{\mathbb{R}}
\newcommand{\Q}{\mathbb{Q}}
\newcommand{\F}{\mathbb{F}}
\newcommand{\Hy}{\mathbb{H}}
\newcommand{\Nov}[2]{\widehat{\Z #1}_{#2}}

\newcommand{\Wh}{\operatorname{Wh}}
\newcommand{\Isom}{\operatorname{Isom}}
\newcommand{\vol}{\operatorname{vol}}
\newcommand{\Mbar}{\overline{M}}
\newcommand{\one}{\mathbf{1}}
\newcommand{\state}{\mathfrak{s}}
\newcommand{\eps}{\varepsilon_1}
\newcommand{\Ip}{I_{\mathrm p}}
\newcommand{\If}{I_{\mathrm f}}
\newcommand{\cG}{\mathcal{G}}
\newcommand{\cT}{\mathcal{T}}

\title[Fibered hyperbolic 7-manifolds in infinitely many commensurability classes]{Infinitely many commensurability classes of hyperbolic \texorpdfstring{$7$}{7}-manifolds that fiber over the circle}
\author{Lizi Guo}
\address{Department of Mathematics, University of Wisconsin--Madison, Madison, WI 53706, USA}
\email{llguo@wisc.edu}

\begin{document}

\begin{abstract}
We construct orientable, finite-volume hyperbolic $7$-manifolds $M_m$ ($m=0$ and all odd $m\ge1$) that fiber smoothly over the circle and are pairwise incommensurable. Each $M_m$ is tessellated by $256$ copies of a right-angled polytope formed by gluing $P^7$ to a chain of $m$ copies of a hyperplane section $Q$ of $P^8$. The fibrations are obtained via Farrell's fibering theorem from Bestvina--Brady Morse functions with simply connected links, proving the vanishing of Novikov homology block by block using exact calculations over the integral group ring, Poincar\'e--Lefschetz duality, and homological perturbation.
\end{abstract}

\maketitle

\section{Introduction}

A smooth manifold \emph{fibers over the circle} if it is the total space of a smooth fiber bundle over $S^1$. By the Gauss--Bonnet theorem, every finite-volume hyperbolic manifold of even dimension has nonzero Euler characteristic and cannot fiber over $S^1$, so the question is of interest in odd dimensions. In dimension $3$, Agol and Wise proved that every finite-volume hyperbolic $3$-manifold has a finite cover fibering over $S^1$ \cite{Agol13,Wise21}. No higher-dimensional example was known until Italiano, Martelli and Migliorini constructed a cusped hyperbolic $5$-manifold fibering over $S^1$ \cite{IMM23}, and it has been conjectured that every odd-dimensional finite-volume hyperbolic manifold is virtually fibered \cite[Conjecture~4.8]{Kielak25}. In dimension $7$, \emph{algebraic} fibrations, namely epimorphisms $\pi_1M\to\Z$ with finitely generated or finitely presented kernel, were found in \cite{IMM24,Fisher23,IMN26}. The existence of a fibered finite-volume hyperbolic $7$-manifold was posed as an open problem in \cite[Section~4]{Kielak25}.

In this paper we construct such manifolds in infinitely many commensurability classes. Two hyperbolic manifolds are \emph{commensurable} if they share a common finite isometric cover.

\begin{maintheorem}\label{thm:main}
Let $m=0$, or let $m\ge1$ be an odd integer. There exists an orientable, finite-volume hyperbolic $7$-manifold $M_m$ with the following properties.
\begin{enumerate}
\item $M_m$ is tessellated by $256$ copies of a right-angled polytope $S_m$ with $56+153m$ facets, where $S_0=P^7$ and for $m\ge1$ the polytope $S_m$ is the union of $P^7$ and $m$ copies of a right-angled polytope $Q$, so that $\vol(M_m)=256\,(\vol(P^7)+m\vol(Q))$. Moreover, $M_m$ has $252+1458m$ cusps, each with a flat $6$-torus cross-section.
\item $M_m$ is the total space of a smooth fiber bundle $M_m\to S^1$.
\item The fiber is diffeomorphic to the interior of a compact, connected, orientable, aspherical $6$-manifold whose boundary is a disjoint union of copies of $T^5$.
\item If $m\ne m'$, then $M_m$ and $M_{m'}$ are not commensurable.
\end{enumerate}
\end{maintheorem}

The case $m=0$ gives an arithmetic example, while $m\ge1$ produces infinitely many non-arithmetic commensurability classes.

\begin{maincorollary}\label{cor:single}
There is an orientable, finite-volume arithmetic hyperbolic $7$-manifold $M$ tessellated by $256$ copies of $P^7$, with $252$ cusps of cross-section $T^6$, that fibers smoothly over $S^1$. Its fiber is diffeomorphic to the interior of a compact, connected, orientable, aspherical $6$-manifold with boundary a disjoint union of copies of $T^5$.
\end{maincorollary}

\begin{maincorollary}\label{cor:classes}
There exist infinitely many commensurability classes of finite-volume hyperbolic $7$-manifolds that fiber over the circle.
\end{maincorollary}

\subsection*{Outline of the proof}
We first construct the manifolds $M_m$ (Section~\ref{sec:manifold}). Starting from the right-angled polytope $P^8\subset\Hy^8$ of \cite{PV05,ERT12}, we take its intersection $Q\subset\Hy^7$ with a symmetry hyperplane. The polytope $Q$ has two ultraparallel facets isometric to $P^6$, matching a facet of $P^7$ \cite{ALR01}. Gluing $P^7$ to a chain of $m$ reflected copies of $Q$ along these facets yields a right-angled polytope $S_m\subset\Hy^7$ whose outward unit normals lie in $\Z[\sqrt2]^8$. Reduction modulo $\sqrt2$ colors the facets of $S_m$ in $V=\F_2^8$, defining an index-$256$ torsion-free normal subgroup $\Gamma_m<\Isom^+(\Hy^7)$ of the reflection group of $S_m$ and the manifold $M_m=\Hy^7/\Gamma_m$ (following \cite{JNW21,IMM23,IMM24}). Since the reflection groups of $P^7$ and $Q$ are incommensurable arithmetic lattices, $M_m$ ($m\ge1$) is a Gromov--Piatetski-Shapiro hybrid \cite{GPS88}, and pairwise incommensurability follows via \cite[Lemma~3.5]{GL14}.

To fiber $M_m$ over the circle, we depart from the $5$-dimensional construction of \cite{IMM23}, which relies on a Bestvina--Brady Morse function whose ascending and descending links are all collapsible so that no critical points occur. In dimension $7$ we do not require collapsible links. Because the dimension is at least $6$, Farrell's fibering theorem \cite{Farrell71,Farrell72}, together with the finite-domination criterion of Ranicki and Sch\"utz \cite{Ranicki95,Schuetz06} and the vanishing of Whitehead groups of torsion-free $\mathrm{CAT}(0)$ groups \cite{BL12}, reduces the fibering of an epimorphism $\varphi\colon\Gamma_m\to\Z$ to three conditions: $\ker\varphi$ is finitely presented, $\varphi$ is nontrivial on each cusp group, and the Novikov homology with coefficients in $\Nov{\Gamma_m}{\pm\varphi}$ vanishes.

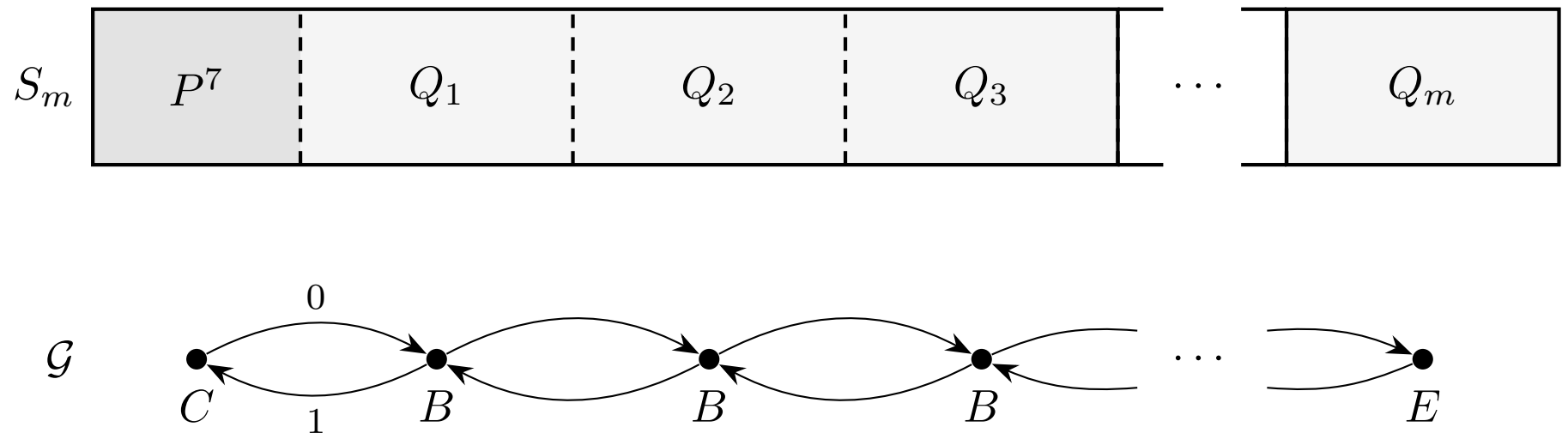
\begin{figure}[ht]
\centering
\begin{tikzpicture}[font=\small,x=1cm,y=1cm]
\fill[gray!22] (0,0) rectangle (1.6,1.2);
\foreach \k in {0,1,2} {\fill[gray!8] (1.6+2.1*\k,0) rectangle (3.7+2.1*\k,1.2);}
\fill[gray!8] (9.2,0) rectangle (11.3,1.2);
\draw[thick] (0,0) rectangle (7.9,1.2); \draw[thick] (9.2,0) rectangle (11.3,1.2);
\draw[thick] (7.9,0)--(8.25,0) (7.9,1.2)--(8.25,1.2) (8.85,0)--(9.2,0) (8.85,1.2)--(9.2,1.2);
\node at (8.55,0.6) {$\cdots$};
\foreach \x in {1.6,3.7,5.8,7.9,9.2} {\draw[thick,densely dashed] (\x,0)--(\x,1.2);}
\node at (0.8,0.6) {$P^7$}; \node at (2.65,0.6) {$Q_1$}; \node at (4.75,0.6) {$Q_2$}; \node at (6.85,0.6) {$Q_3$}; \node at (10.25,0.6) {$Q_m$};
\node[left] at (0,0.6) {$S_m$};
\foreach \x/\n in {0.8/0,2.65/1,4.75/2,6.85/3,10.25/m} {\node[circle,fill=black,inner sep=1.6pt] (v\n) at (\x,-1.5) {};}
\node[below=3pt] at (v0) {$C$}; \node[below=3pt] at (v1) {$B$}; \node[below=3pt] at (v2) {$B$}; \node[below=3pt] at (v3) {$B$}; \node[below=3pt] at (vm) {$E$};
\foreach \a/\b in {v0/v1,v1/v2,v2/v3}
 {\draw[-{Stealth[length=2mm]}] (\a) to[bend left=28] (\b); \draw[-{Stealth[length=2mm]}] (\b) to[bend left=28] (\a);}
\draw (v3) to[bend left=14] ($(v3)+(1.2,0.22)$); \draw[-{Stealth[length=2mm]}] ($(v3)+(1.2,-0.22)$) to[bend left=14] (v3);
\draw[-{Stealth[length=2mm]}] ($(vm)+(-1.2,0.22)$) to[bend left=14] (vm); \draw (vm) to[bend left=14] ($(vm)+(-1.2,-0.22)$);
\node at (8.55,-1.5) {$\cdots$};
\node[left] at (0,-1.5) {$\cG$};
\node[font=\scriptsize] at (1.72,-1.02) {$0$}; \node[font=\scriptsize] at (1.72,-1.98) {$1$};
\end{tikzpicture}
\caption{Top: the polytope $S_m$ is the union of $P^7$ and $m$ copies of $Q$, glued along facets isometric to $P^6$ (dashed). Bottom: the graph $\cG$ of the block decomposition of $M_m$. Each dashed facet gives a hypersurface with two components, labelled $0$ and $1$. An edge points to the block for which the component is incoming.}\label{fig:chain}
\end{figure}

In Section~\ref{sec:character} we choose a state and positive integer weights on the facets of $S_m$. These determine a character $\alpha\colon\Gamma_m\to\Z$ and an $\alpha$-equivariant Morse function on the Davis complex. For both $P^7$ and $Q$, every ascending and descending link is simply connected, and its reduced homology lives only in degrees $2$ and $3$, with total rank $24$ for $P^7$ and $428$ for $Q$, restricting to $214$ on the index-$2$ subspace $U\subset V$ (Proposition~\ref{prop:links}). Rather than dealing with a single matrix whose size grows with $m$, we cut $M_m$ along totally geodesic hypersurfaces into $m+1$ \emph{blocks}, namely a cap $C$ built from $P^7$, $m-1$ isometric bridges $B$, and an end block $E$ built from $Q$ (Figure~\ref{fig:chain}). Each cutting hypersurface has two components, and we orient the dual graph $\cG$ so that each component is incoming for exactly one of the two blocks it bounds.

Section~\ref{sec:novikov} establishes that each block has vanishing Novikov homology relative to its incoming boundary. Contracting the vertex complexes locally and applying homological perturbation reduces the question to the invertibility of a $24\times24$ matrix for the cap $C$, and of a $214\times214$ matrix for a half-sized core $X$ inside $Q$. After a $29$-column unitriangular change of basis over $\Z\Gamma_X$, both matrices are inverted by a leading-term dominance argument (Proposition~\ref{lem:dominance}). Relative acyclicity of $M_0$, $B$, and $E$ then follows from $C$ and $X$ by Poincar\'e--Lefschetz duality and by restoring the cut facets with sufficiently large weights (Lemma~\ref{lem:restore}).

In Section~\ref{sec:gluing} we glue the blocks back together. The edge orientations on $\cG$ lift to the Bass--Serre tree $\cT$ and produce a character $\psi\colon\Gamma_m\to\Z$ vanishing on all block and hypersurface stabilizers. When $N$ is sufficiently large, the perturbed character $\alpha_N=\alpha+N\psi$ is surjective, nontrivial on every cusp, has finitely presented kernel, and has vanishing Novikov homology for both $\pm\alpha_N$ (Theorem~\ref{thm:glue}). Farrell's theorem applied to the compact core $\Mbar\subset M_m$ then produces the fibration. The restriction to odd $m\ge1$ is for convenience (Remark~\ref{rem:even}).

\section{Hybrid right-angled polytopes and commensurability}\label{sec:manifold}

Let $\R^{7,1}$ denote $\R^8$ with coordinates $x=(x_1,\dots,x_7;t)$ and Lorentzian form $\langle x,y\rangle=\sum_{k=1}^7x_ky_k-tt'$, and let $\Hy^7=\{x\in\R^{7,1}:\langle x,x\rangle=-1,\ t>0\}$. A unit vector $e$ ($\langle e,e\rangle=1$) defines the half-space $\{x\in\Hy^7:\langle x,e\rangle\le0\}$ with outward unit normal $e$. Two half-spaces with normals $e\ne\pm e'$ have orthogonal boundaries if $\langle e,e'\rangle=0$, boundaries tangent at infinity if $\langle e,e'\rangle=\pm1$, and \emph{ultraparallel} boundaries if $|\langle e,e'\rangle|>1$.

\subsection{\texorpdfstring{The polytopes $P^7$ and $Q$}{The polytopes P7 and Q}}
Let $e_0,\dots,e_{55}\in\Z^8$ be the vectors of Table~\ref{tab:normals}, obtained from
\[
(-1,0,0,0,0,0,0;0),\quad(1,1,0,0,0,0,0;1),\quad(1,1,1,1,1,0,0;2),\quad(2,1,1,1,1,1,1;3)
\]
by permuting the first seven coordinates and negating $x_1$, and define
\[
P^7=\{x\in\Hy^7:\ \langle x,e_j\rangle\le0\ \text{for}\ j=0,\dots,55\}.
\]
Up to isometry, this is the right-angled polytope of Agol, Long and Reid \cite{ALR01}, see also \cite{PV05,ERT12}. We number its facets $F_0,\dots,F_{55}$ as in Table~\ref{tab:normals}, so $F_8$ has outward normal $e_8=(1,0,\dots,0;0)$. Whenever $i\ne j$ we have $\langle e_i,e_j\rangle\in\{0,-1,-2\}$, and we write $i\sim j$ when $\langle e_i,e_j\rangle=0$. We recall the basic facts about $P^7$ from \cite{ALR01,ERT12,IMM24} and \cite{PV05}.

\begin{lemma}\label{lem:P7}
$P^7$ is a finite-volume right-angled polytope with $56$ facets, each meeting $27$ others. Two facets $F_i$ and $F_j$ intersect if and only if $i\sim j$, and then they meet orthogonally. $P^7$ has $576$ finite vertices and $126$ ideal vertices, represented by the null vectors $e_i+e_j$ with $\langle e_i,e_j\rangle=-1$. Each ideal vertex lies on $12$ facets forming $6$ tangent pairs, its link is a Euclidean $6$-cube, and $F_8=P^7\cap\{x_1=0\}$ is isometric to $P^6$.
\end{lemma}

In $\R^{8,1}$, let $\mathcal N_8\subset\Z^9$ be the set of $240$ vectors obtained by permuting the first eight coordinates of
\begin{gather*}
(-1,0,0,0,0,0,0,0;0),\quad(1,1,0,0,0,0,0,0;1),\quad(1,1,1,1,1,0,0,0;2),\quad(2,1,1,1,1,1,1,0;3),\\
(2,2,2,1,1,1,1,1;4),\quad(2,2,2,2,2,2,1,1;5),\quad(3,2,2,2,2,2,2,2;6),
\end{gather*}
defining the finite-volume right-angled polytope $P^8=\{x\in\Hy^8:\langle x,n\rangle\le0\ \text{for all}\ n\in\mathcal N_8\}$ \cite{PV05,ERT12}. The involution $\sigma$ swapping $x_1$ and $x_2$ preserves $\mathcal N_8$, and its fixed hyperplane $\{x_1=x_2\}$ is identified isometrically with $\R^{7,1}$ via
\[
(x,x,x_3,\dots,x_8;t)\longmapsto(\sqrt2\,x,x_3,\dots,x_8;t).
\]
Set $Q=P^8\cap\{x_1=x_2\}\subset\Hy^7$. Since $(1,\dots,1;3)\in\{x_1=x_2\}$ satisfies all defining inequalities of $P^8$ strictly, $Q$ has nonempty interior in $\Hy^7$. Its ideal boundary is contained in the finite set of ideal vertices of $P^8$, so $Q$ is the convex hull of finitely many finite and ideal points and therefore has finite volume in $\Hy^7$.

For $n\in\mathcal N_8$ and $x\in\{x_1=x_2\}$, we have $\langle x,n\rangle=\langle x,\frac12(n+\sigma n)\rangle$. There are three cases:
\begin{itemize}
\item $\sigma n=n$ ($126$ vectors): $n=(a,a,b_3,\dots,b_8;c)$ yields the unit normal $(\sqrt2\,a,b_3,\dots,b_8;c)\in\R^{7,1}$.
\item $\sigma n\ne n$ and $\langle n,\sigma n\rangle=0$ ($56$ pairs): $n=(a,b,b_3,\dots,b_8;c)$ yields the unit normal $(a+b,\sqrt2\,b_3,\dots,\sqrt2\,b_8;\sqrt2\,c)\in\R^{7,1}$.
\item $\sigma n\ne n$ and $\langle n,\sigma n\rangle\ne0$: only $n=(2,0,1,1,1,1,1,1;3)$ and $\sigma n$, for which $\frac12(n+\sigma n)=(1,\dots,1;3)$ is timelike and the inequality is redundant on $\Hy^8$.
\end{itemize}
We order the resulting $182$ unit vectors $g_0,\dots,g_{181}\in\Z[\sqrt2]^8$ lexicographically, where $p+q\sqrt2$ precedes $p'+q'\sqrt2$ if $p<p'$, or if $p=p'$ and $q<q'$. Then $Q=\{x\in\Hy^7:\langle x,g_i\rangle\le0\ \text{for}\ 0\le i\le181\}$, with facets $E_i=Q\cap g_i^\perp$. The vectors of the first kind are indexed by $\If=\{1,\dots,126\}$ (\emph{fixed} facets) and those of the second kind by $\Ip=\{0,127,\dots,181\}$ (\emph{paired} facets). The extreme normals are
\[
g_0=(-1,0,0,0,0,0,0;0),\qquad g_{181}=(5,2\sqrt2,2\sqrt2,2\sqrt2,2\sqrt2,2\sqrt2,2\sqrt2;6\sqrt2).
\]
Again we write $i\sim j$ if $i\ne j$ and $\langle g_i,g_j\rangle=0$.

\begin{lemma}\label{lem:Q}
\begin{enumerate}
\item $Q$ is a finite-volume right-angled polytope with $182$ facets $E_0,\dots,E_{181}$. For $i\ne j$, the inner product $\langle g_i,g_j\rangle$ lies in $\{0,-1,-\sqrt2,-2,-2\sqrt2,-3,-5\}$, and $E_i\cap E_j\ne\emptyset$ if and only if $i\sim j$, in which case the facets are orthogonal.
\item Each fixed facet meets $44$ other facets. Each paired facet meets $27$ fixed facets and no paired facet.
\item The numbers of $k$-tuples of pairwise intersecting facets are $182$, $3528$, $22176$, $60480$, $76608$, $40320$, $4608$ for $k=1,\dots,7$, and $0$ for $k=8$. In particular, $Q$ has $4608$ finite vertices.
\item $Q$ has $756$ ideal vertices, represented by the null vectors $g_i+g_j$ with $\langle g_i,g_j\rangle=-1$. Each ideal vertex lies on $12$ facets that form $6$ tangent pairs, namely one pair of paired facets and five pairs of fixed facets.
\item We have $\langle g_0,g_{181}\rangle=-5$, so $E_0$ and $E_{181}$ are ultraparallel and no facet meets both. The facets meeting $E_0$ are $E_1,\dots,E_{27}$, and their normals are the vectors $e_j$ with $j\sim8$ in Table~\ref{tab:normals}. Hence $E_0=Q\cap\{x_1=0\}$ coincides with $F_8=P^7\cap\{x_1=0\}$.
\item The central reflection $\jmath(x)=-x-2\langle c,x\rangle c$ at $c=(\sqrt2,1,1,1,1,1,1;3)\in\Hy^7$ preserves $Q$, preserves $\If$ and $\Ip$, and sends $g_0$ to $g_{181}$.
\end{enumerate}
\end{lemma}

\begin{proof}
The inner products, incidences, and the action of $\jmath$ on $\{g_i\}$ are verified directly (Appendix~\ref{app:computations}). Since $\langle g_i,g_j\rangle$ is either $0$ or at most $-1$ for $i\ne j$, the polytope $Q$ is right-angled. By \cite[Theorem~3.1]{Vinberg85}, its codimension-$k$ faces correspond to $k$-tuples of pairwise orthogonal normals, which proves (1)--(3). The ideal vertices are the lightlike rays in $\partial\Hy^7$ that satisfy all defining inequalities and are orthogonal to a tangent pair ($\langle g_i,g_j\rangle=-1$). This gives (4), while (5) and (6) follow directly from the computation.
\end{proof}

We call $E_0,E_{181}$ of $Q$ and $F_8$ of $P^7$ \emph{ports}.

\subsection{\texorpdfstring{The hybrid chain $S_m$ and its $256$-sheeted cover $M_m$}{The hybrid chain Sm and its 256-sheeted cover Mm}}
By Lemma~\ref{lem:Q}(5), $P^7\subset\{x_1\le0\}$ and $Q\subset\{x_1\ge0\}$ meet precisely along $F_8=E_0$. Let $r_i$ denote the reflection in $E_i$, and put
\[
h_k=\underbrace{r_{181}\,r_0\,r_{181}\cdots}_{k-1\ \text{factors}},\qquad Q_k=h_k(Q)\qquad(k\ge1).
\]
Then $Q_1=Q$, and $Q_{k+1}$ is the reflection of $Q_k$ in its \emph{outlet} $h_k(E_o)$, where $o=181$ if $k$ is odd and $o=0$ if $k$ is even. The opposite port of $Q_k$ is its \emph{inlet}. Set $S_0=P^7$ and $S_m=P^7\cup Q_1\cup\dots\cup Q_m$ for $m\ge1$.

\begin{lemma}\label{lem:Sm}
For $m\ge1$, $S_m$ is a finite-volume right-angled polytope with $56+153m$ facets and $126+729m$ ideal vertices, all of whose outward unit normals lie in $\Z[\sqrt2]^8$. Its facets consist of the unglued facets of $P^7,Q_1,\dots,Q_m$, with each pair of facets from consecutive pieces meeting their common port merging into a single facet of $S_m$ (called \emph{parts} of that facet).
\end{lemma}

\begin{proof}
Consecutive pieces lie on opposite sides of their common port and meet along it (by Lemma~\ref{lem:Q}(5) for $P^7,Q_1$ and reflection symmetry for $Q_k,Q_{k+1}$). Facets meeting a port are orthogonal to it, so adjacent facets of consecutive pieces lie in common hyperplanes and form a locally convex union. Since the two ports of each $Q_k$ are ultraparallel and share no adjacent facet, induction on $m$ shows that $S_m$ is a right-angled convex polytope with $55+m(182-2-27)+1=56+153m$ facets and $126+m(756-27)=126+729m$ ideal vertices (each port contains $27$ ideal vertices). Since reflections in unit vectors of $\Z[\sqrt2]^8$ preserve $\Z[\sqrt2]^8$, all outward normals lie in $\Z[\sqrt2]^8$.
\end{proof}

Let $V=\F_2^8$ with standard bilinear form $u\cdot v=\sum_k u_kv_k$, and let $\Z[\sqrt2]^8\to V$ be coordinatewise reduction modulo $\sqrt2$ ($p+q\sqrt2\mapsto p\bmod2$), which intertwines $\langle\,,\rangle$ and $\cdot$. Set $\eps=(1,0,\dots,0)$, $\one=(1,\dots,1)$, and $U=\{v\in V:v\cdot\eps=0\}$. For any finite-volume right-angled polytope $P\subset\Hy^7$ with facets $F_i$ ($i\in I$) and outward unit normals $e_i\in\Z[\sqrt2]^8$, define the \emph{color} $c_i\in V$ of $F_i$ as the reduction of $e_i$ modulo $\sqrt2$. Then
\begin{equation}\label{eq:colors}
c_i\cdot c_i=c_i\cdot\one=1,\qquad c_i\cdot c_j=0\quad\text{whenever }F_i\cap F_j\ne\emptyset.
\end{equation}
For $Q$, paired facets have color $\eps$ while fixed facets have colors in $U$. Because reflections in unit vectors of $\Z[\sqrt2]^8$ reduce to the identity modulo $\sqrt2$, $h_k(E_i)$ has the same color as $E_i$, and merged facets of consecutive pieces in $S_m$ carry equal colors.

\begin{figure}[ht]
\centering
\begin{tikzpicture}[font=\small]
\fill[gray!12] (-2.1,-1.5) rectangle (2.1,1.5);
\draw[thick] (-2.1,0)--(2.1,0); \draw[thick] (0,-1.5)--(0,1.5);
\node at (-1.05,-0.75) {$v$};
\node at (1.05,-0.75) {$v+c_i$};
\node at (-1.05,0.75) {$v+c_j$};
\node at (1.05,0.75) {$v+c_i+c_j$};
\node[above] at (0,1.5) {$F_i$}; \node[right] at (2.1,0) {$F_j$};
\fill (0,0) circle (1.6pt);
\end{tikzpicture}
\caption{The labels of the four copies of the polytope around a codimension-two face $F_i\cap F_j$. Crossing $F_i$ adds $c_i$ to the label, and crossing $F_j$ adds $c_j$.}\label{fig:moves}
\end{figure}
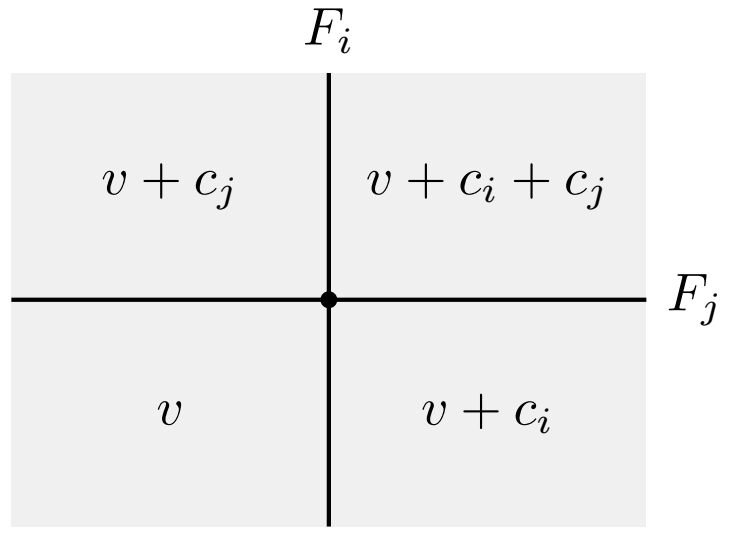

Let $W$ be the right-angled Coxeter group generated by the facet reflections $r_i$ ($i\in I$), with nerve $L$ (the flag complex on $I$ whose simplices, called \emph{cliques}, are sets of pairwise intersecting facets). By Poincar\'e's polyhedron theorem \cite{Vinberg85}, $W$ acts on $\Hy^7$ with fundamental domain $P$. Define $\rho\colon W\to V$ by $\rho(r_i)=c_i$. For $P=S_m$, write $W_m,\rho_m$ and put
\[
\Gamma_m=\ker\rho_m,\qquad M_m=\Hy^7/\Gamma_m .
\]
Labeling each tile $wS_m$ by $\rho_m(w)\in V$, $M_m$ is assembled from $256$ copies of $S_m$ indexed by $V$, gluing copy $v$ to copy $v+c_i$ along $F_i$ (Figure~\ref{fig:moves}).

\begin{proposition}\label{prop:cover}
\begin{enumerate}
\item $\Gamma_m$ is a torsion-free normal subgroup of index $256$ in $W_m$ contained in $\Isom^+(\Hy^7)$.
\item $M_m$ is an orientable, finite-volume hyperbolic $7$-manifold tessellated by $256$ copies of $S_m$, with $\vol(M_m)=256\,(\vol(P^7)+m\vol(Q))$.
\item $M_m$ has $252+1458m$ cusps, each with a flat $6$-torus cross-section.
\end{enumerate}
\end{proposition}

\begin{proof}
(1) The colors of the $55$ facets of $P^7\setminus\{F_8\}$ span $V$ (Appendix~\ref{app:computations}), so $\rho_m$ is surjective and $[W_m:\Gamma_m]=256$. By \eqref{eq:colors}, the colors of any clique $S$ are orthonormal, hence linearly independent, so $\rho_m$ is injective on $W_S\cong(\Z/2)^{|S|}$. Every finite subgroup of $W_m$ is conjugate into some $W_S$ \cite[Chapter~12]{Davis08}. Since $V$ is abelian, $\Gamma_m$ is torsion-free. The composition of $\rho_m$ with $v\mapsto v\cdot\one$ is the orientation character $r_i\mapsto1$, so $\Gamma_m<\Isom^+(\Hy^7)$. Part~(2) follows immediately.

(3) An ideal vertex $\xi$ of $S_m$ has Euclidean $6$-cube link, so it lies on six tangent pairs $(F_{\alpha_k},F_{\beta_k})$ with facets from different pairs orthogonal. Its stabilizer is $W_\xi=\prod_{k=1}^6\langle r_{\alpha_k},r_{\beta_k}\rangle\cong D_\infty^6$, with translation subgroup $\Z^6=\langle t_1,\dots,t_6\rangle$, where $t_k=r_{\alpha_k}r_{\beta_k}$. Put $\delta_k=c_{\alpha_k}+c_{\beta_k}$. Since $c_{\alpha_k}\cdot c_{\beta_k}=1$, the vectors $c_{\alpha_1},\dots,c_{\alpha_6}$ are orthonormal, and the $\delta_k$ lie in their orthogonal complement $W^\perp$ with $\delta_k\cdot\delta_l=0$. For $g=\prod_k t_k^{z_k}r_{\alpha_k}^{\epsilon_k}\in W_\xi$, pairing $\rho_m(g)=\sum_k z_k\delta_k+\sum_k\epsilon_k c_{\alpha_k}$ with $c_{\alpha_l}$ gives $\epsilon_l$. Thus $\Gamma_m\cap W_\xi$ is a finite-index sublattice of $\Z^6$. Because $\Gamma_m$ is normal in $W_m$, the same holds for every conjugate of $W_\xi$.

The orthogonal complement $W^\perp$ is a nondegenerate $2$-plane in $V$, so its totally isotropic subspace $\operatorname{span}\{\delta_k\}$ has dimension at most $1$. At least five of the six tangent pairs at $\xi$ come from non-port facets of a single piece ($P^7$ or some $Q_k$). Every tangent pair of $P^7$ and every fixed tangent pair of $Q$ satisfies $c_i+c_j\ne0$ (Appendix~\ref{app:computations}). Hence $|\rho_m(W_\xi)|=2^7$, and there are $[V:\rho_m(W_\xi)]=2$ cusps over each ideal vertex. By Lemmas~\ref{lem:P7} and~\ref{lem:Sm}, the total number of cusps is $2(126+729m)=252+1458m$.
\end{proof}

\subsection{Incommensurability}
Let $p\colon\Hy^7\to M_m$ be the covering projection, and set $A_m=p(W_m\cdot P^7)$ and $T_m=p(W_m\cdot(Q_1\cup\dots\cup Q_m))$ (with $T_0=\emptyset$). Then $M_m=A_m\cup T_m$ along their common totally geodesic boundary, with
\begin{equation}\label{eq:volumes}
\vol(A_m)=256\vol(P^7),\qquad\vol(T_m)=256\,m\vol(Q).
\end{equation}

\begin{lemma}\label{lem:arith}
\begin{enumerate}
\item The reflection groups of $P^7$ and $Q$ are commensurable with $O(q_P,\Z)$ and $O(q_Q,\Z)$, respectively, where
\[
q_P=x_1^2+x_2^2+\dots+x_7^2-t^2,\qquad q_Q=2x_1^2+x_2^2+\dots+x_7^2-t^2 .
\]
These two arithmetic groups are not commensurable in $\Isom(\Hy^7)$, even up to conjugacy.
\item $A_m$ (resp.\ $T_m$) embeds isometrically as a submanifold with totally geodesic boundary in $\Hy^7/\Lambda$ for a torsion-free finite-index subgroup $\Lambda$ of the reflection group of $P^7$ (resp.\ $Q$).
\end{enumerate}
\end{lemma}

\begin{proof}
(1) The reflections in $e_j\in\Z^8$ preserve $\Z^8$ and $q_P$. In the coordinates $(x_1/\sqrt2,x_2,\dots,x_7;t)$, the form $\langle\,,\rangle$ becomes $q_Q$. Each $g_i$ ($i\in\If$) is then integral of norm $1$, and each $g_i$ ($i\in\Ip$) is $w/\sqrt2$ for an integral vector $w$ with $q_Q(w)=2$ whose pairings with all integral vectors are even. Therefore all reflections in $g_i$ preserve $\Z^8$ and $q_Q$. Being lattices, both reflection groups have finite index in $O(q_P,\Z)$ and $O(q_Q,\Z)$. By \cite[2.6]{GPS88}, they are commensurable up to conjugacy if and only if $q_P$ and $q_Q$ are similar over $\Q$. Since both forms have $8$ variables, the discriminant in $\Q^\times/(\Q^\times)^2$ is invariant under similarity, and it is $-1$ for $q_P$ versus $-2$ for $q_Q$.

(2) If $m=0$ we take $\Lambda=\Gamma_0$. For $m\ge1$, let $R$ be $P^7$ or $Q_1\cup\dots\cup Q_m$, with port $F=F_8$ and reflection group $W_R$. Define $\rho'\colon W_R\to V\oplus\F_2$ by $r_i\mapsto(c_i,0)$ for $i\ne F$ and $r_F\mapsto(0,1)$. As in Proposition~\ref{prop:cover}(1), $\Lambda=\ker\rho'$ is torsion-free of finite index, and $\Hy^7/\Lambda$ is the double of $A_m$ (resp.\ $T_m$) along its boundary.
\end{proof}

\begin{theorem}\label{prop:commens}
If $M_m$ and $M_{m'}$ are commensurable, then $m=m'$.
\end{theorem}

\begin{proof}
Let $\pi\colon N\to M_m$ and $\pi'\colon N\to M_{m'}$ be finite covers of degrees $d$ and $d'$. By Lemma~\ref{lem:arith} and \cite[Lemma~3.5]{GL14}, which holds in dimension at least $4$, the sets $\pi^{-1}(A_m)$ and $\pi'^{-1}(T_{m'})$ meet in a set with empty interior. Any point of this intersection lies on the boundary of one of the two pieces, so the intersection is null. Thus $\pi^{-1}(A_m)$ is contained in $\pi'^{-1}(A_{m'})$ up to a null set, and by symmetry the two have equal volume. Using \eqref{eq:volumes} gives $d=d'$, and comparing $\vol(N)$ yields $m=m'$.
\end{proof}

\section{Bestvina--Brady Morse theory, blocks, and link topology}\label{sec:character}

In this section we work with a finite-volume right-angled polytope $P\subset\Hy^7$ as in Section~\ref{sec:manifold}. Write $F_i$ ($i\in I$) for its facets, $c_i\in V$ for their colors, $W$ for the reflection group, $L$ for its nerve, $\rho\colon W\to V$ for the coloring homomorphism, and $\Gamma=\ker\rho$.

\subsection{The Davis complex and height functions}
Dual to the tessellation $\{wP\}$ is the \emph{Davis complex} $\Sigma$ of $W$, the cube complex with vertex set $W$ and cubes $wW_S$ for cliques $S$ of $L$. Each vertex link is a copy of the flag complex $L$, so $\Sigma$ is a $\mathrm{CAT}(0)$ cube complex and in particular contractible \cite[Chapter~12]{Davis08}. Because $\Gamma$ acts freely and cocompactly on $\Sigma$, the quotient $\Sigma/\Gamma$ is a finite $K(\Gamma,1)$.

Fix a \emph{state} $\state\in\F_2^I$ and positive integer \emph{weights} $a_i$ ($i\in I$). For each $v\in V$, define
\[
D(v)=\{i\in I:\ \state_i+v\cdot c_i=1\},\qquad
\alpha(v,i)=\begin{cases}-a_i,& i\in D(v),\\ \phantom{-}a_i,& i\notin D(v).\end{cases}
\]

\begin{lemma}\label{lem:cocycle}
There is a unique map $f\colon W\to\Z$ with $f(1)=0$ and $f(wr_i)-f(w)=\alpha(\rho(w),i)$ for all $w\in W$ and $i\in I$. Moreover, there is a homomorphism $\varphi\colon\Gamma\to\Z$ with $f(\gamma w)=\varphi(\gamma)+f(w)$ for $\gamma\in\Gamma$ and $w\in W$. If $\tau\in W$ satisfies $\rho(\tau)=\one$, then $\varphi(\tau\gamma\tau^{-1})=-\varphi(\gamma)$ for all $\gamma\in\Gamma$.
\end{lemma}

\begin{proof}
By \eqref{eq:colors}, we have $\alpha(v+c_i,i)=-\alpha(v,i)$, and $\alpha(v+c_i,j)=\alpha(v,j)$ whenever $F_i\cap F_j\ne\emptyset$. The increments therefore sum to $0$ around each relator $r_i^2$ and $(r_ir_j)^2$, so $f$ is well defined. For $\gamma\in\Gamma$ we have $\rho(\gamma w)=\rho(w)$, so $f(\gamma w)-f(w)$ is independent of $w$ and defines $\varphi$. Finally, $c_i\cdot\one=1$ implies $D(v+\one)=I\setminus D(v)$ and $\alpha(v+\one,i)=-\alpha(v,i)$. Thus $f(\tau w)=-f(w)+c$ for a constant $c$, which gives $\varphi(\tau\gamma\tau^{-1})=-\varphi(\gamma)$.
\end{proof}

Recall from \cite[Definitions~2.1--2.4]{BB97} that a \emph{Morse function} on an affine cell complex $Y$ is a map $f\colon Y\to\R$ that is affine and nonconstant on every cell of positive dimension, with discrete image on $Y^{(0)}$. The \emph{descending} (resp.\ \emph{ascending}) \emph{link} of $y\in Y^{(0)}$ is the subcomplex of $\operatorname{lk}(y)$ formed by the cells on which $f$ attains its maximum (resp.\ minimum) at $y$. By \cite[Theorem~4.1(3)]{BB97}, if a group $G$ acts freely, properly, cocompactly and affinely on a contractible affine cell complex $Y$ and $f\colon Y\to\R$ is a Morse function equivariant with respect to an epimorphism $\varphi\colon G\to\Z$ such that all ascending and descending links are simply connected, then $\ker\varphi$ is finitely presented.

\begin{figure}[ht]
\centering
\begin{tikzpicture}[font=\small,vx/.style={circle,fill=black,inner sep=1.4pt},lab/.style={align=center,font=\small,inner sep=2pt}]
\draw[-{Stealth[length=2mm]}] (-3.3,-2.3)--(-3.3,0.3) node[above] {$f$};
\begin{scope}
\coordinate (w) at (0,0); \coordinate (wi) at (-1.2,-1.3); \coordinate (wj) at (1.4,-0.9); \coordinate (wij) at (0.2,-2.2);
\fill[gray!15] (w)--(wi)--(wij)--(wj)--cycle;
\begin{scope}
\clip (w)--(wi)--(wij)--(wj)--cycle;
\foreach \y in {-0.3,-0.6,...,-2.15} {\draw[gray!75,thin] (-2,\y)--(2,\y);}
\end{scope}
\draw[thick] (w)--(wi)--(wij)--(wj)--cycle;
\draw[-{Stealth[length=2mm]},thick] (w)--($(w)!0.6!(wi)$);
\draw[-{Stealth[length=2mm]},thick] (w)--($(w)!0.6!(wj)$);
\draw[-{Stealth[length=2mm]},thick] (wi)--($(wi)!0.6!(wij)$);
\draw[-{Stealth[length=2mm]},thick] (wj)--($(wj)!0.6!(wij)$);
\foreach \p in {w,wi,wj,wij} {\node[vx] at (\p) {};}
\node[lab,above=2pt] at (w) {$w$\\ $\scriptstyle f(w)$};
\node[lab,left=2pt] at (wi) {$wr_i$\\ $\scriptstyle f(w)-a_i$};
\node[lab,right=2pt] at (wj) {$wr_j$\\ $\scriptstyle f(w)-a_j$};
\node[lab,below=2pt] at (wij) {$wr_ir_j$\\ $\scriptstyle f(w)-a_i-a_j$};
\node at (0.1,-3.55) {(a) $i,j\in D(v)$};
\end{scope}
\begin{scope}[shift={(5.6,-1.1)}]
\coordinate (w) at (0,0); \coordinate (wi) at (1.2,-1.3); \coordinate (wj) at (1.4,0.9); \coordinate (wij) at (2.6,-0.4);
\fill[gray!15] (w)--(wi)--(wij)--(wj)--cycle;
\begin{scope}
\clip (w)--(wi)--(wij)--(wj)--cycle;
\foreach \y in {-1.1,-0.8,...,0.8} {\draw[gray!75,thin] (-1,\y)--(4,\y);}
\end{scope}
\draw[thick] (w)--(wi)--(wij)--(wj)--cycle;
\draw[-{Stealth[length=2mm]},thick] (w)--($(w)!0.6!(wi)$);
\draw[-{Stealth[length=2mm]},thick] (wj)--($(wj)!0.6!(w)$);
\draw[-{Stealth[length=2mm]},thick] (wij)--($(wij)!0.6!(wi)$);
\draw[-{Stealth[length=2mm]},thick] (wj)--($(wj)!0.6!(wij)$);
\foreach \p in {w,wi,wj,wij} {\node[vx] at (\p) {};}
\node[lab,left=2pt] at (w) {$w$\\ $\scriptstyle f(w)$};
\node[lab,below=2pt] at (wi) {$wr_i$\\ $\scriptstyle f(w)-a_i$};
\node[lab,above=2pt] at (wj) {$wr_j$\\ $\scriptstyle f(w)+a_j$};
\node[lab,right=2pt] at (wij) {$wr_ir_j$\\ $\scriptstyle f(w)-a_i+a_j$};
\node at (1.3,-2.45) {(b) $i\in D(v)$, $j\notin D(v)$};
\end{scope}
\end{tikzpicture}
\caption{The extension of $f$ to a square $wW_S$ with $S=\{i,j\}$, drawn with the height increasing upwards. The horizontal segments are level sets of $f$, and the arrows point towards the lower endpoint of each edge. Here $v=\rho(w)$. (a) If $i,j\in D(v)$, then $f$ attains its maximum on the square at $w$. (b) If $i\in D(v)$ and $j\notin D(v)$, then the maximum is at $wr_j$ and the minimum is at $wr_i$.}\label{fig:extension}
\end{figure}
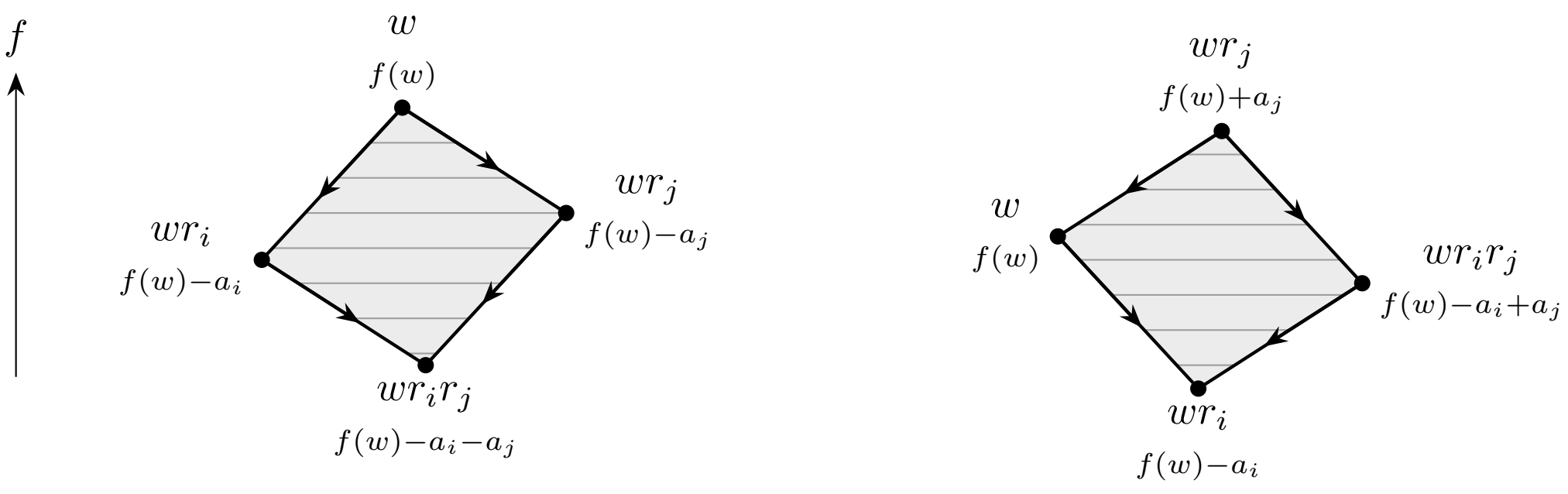

On any cube $wW_S$ of $\Sigma$, distinct $i,j\in S$ satisfy $c_i\cdot c_j=0$, so all edges in direction $j$ have the same increment $\pm a_j$. Hence $f$ extends uniquely to an affine map $f\colon\Sigma\to\R$ (Figure~\ref{fig:extension}). Let $L_v$ be the full subcomplex of $L$ spanned by $D(v)$.

\begin{lemma}\label{lem:morse}
The extension $f\colon\Sigma\to\R$ is a $\varphi$-equivariant Morse function. For $w\in W$ with $v=\rho(w)$, $f$ attains its maximum (resp.\ minimum) on $wW_S$ at $w$ if and only if $S\subseteq D(v)$ (resp.\ $S\cap D(v)=\emptyset$). In particular, the descending and ascending links of $w$ are $L_v$ and $L_{v+\one}$.
\end{lemma}

\begin{proof}
By construction, $f$ is affine on cubes, integer-valued on vertices, and changes by $\pm a_j\ne0$ along each edge, so it is a Morse function. Equivariance extends from the vertices by affinity. An affine function on $wW_S$ is maximized (resp.\ minimized) at $w$ if and only if it decreases (resp.\ increases) along every edge $w\to wr_j$ ($j\in S$), that is, if and only if $S\subseteq D(v)$ (resp.\ $S\subseteq I\setminus D(v)=D(v+\one)$). Under $\operatorname{lk}(w)\cong L$, the descending and ascending links are therefore $L_v$ and $L_{v+\one}$.
\end{proof}

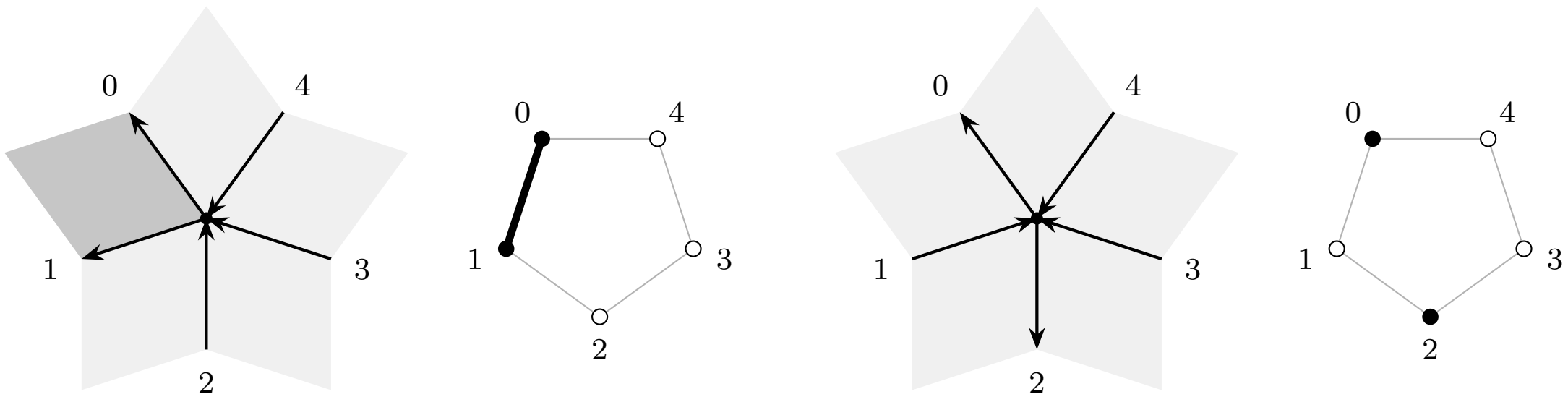
\begin{figure}[ht]
\centering
\begin{tikzpicture}[vx/.style={circle,fill=black,inner sep=1.5pt},ov/.style={circle,draw,fill=white,inner sep=1.4pt}]
\begin{scope}
\foreach \k in {0,...,4} {\coordinate (e\k) at ({126+72*\k}:1.2);}
\foreach \k [evaluate=\k as \l using {int(mod(\k+1,5))}] in {0,...,4} {\fill[gray!12] (0,0)--(e\k)--($(e\k)+(e\l)$)--(e\l)--cycle;}
\fill[gray!45] (0,0)--(e0)--($(e0)+(e1)$)--(e1)--cycle;
\foreach \k in {0,1} {\draw[-{Stealth[length=2mm]},thick] (0,0)--(e\k);}
\foreach \k in {2,3,4} {\draw[{Stealth[length=2mm]}-,thick] (0,0)--(e\k);}
\foreach \k in {0,...,4} {\node[font=\scriptsize] at ({126+72*\k}:1.5) {$\k$};}
\fill (0,0) circle (1.6pt);
\end{scope}
\begin{scope}[xshift=3.6cm]
\foreach \k in {0,...,4} {\coordinate (n\k) at ({126+72*\k}:0.9);}
\draw[gray!60] (n0)--(n1)--(n2)--(n3)--(n4)--cycle;
\draw[line width=2pt] (n0)--(n1);
\foreach \k in {0,1} {\node[vx] at (n\k) {};}
\foreach \k in {2,3,4} {\node[ov] at (n\k) {};}
\foreach \k in {0,...,4} {\node[font=\scriptsize] at ({126+72*\k}:1.2) {$\k$};}
\end{scope}
\node[font=\small] at (1.8,-2.2) {(a) $D(v)=\{0,1\}$};
\begin{scope}[xshift=7.6cm]
\foreach \k in {0,...,4} {\coordinate (e\k) at ({126+72*\k}:1.2);}
\foreach \k [evaluate=\k as \l using {int(mod(\k+1,5))}] in {0,...,4} {\fill[gray!12] (0,0)--(e\k)--($(e\k)+(e\l)$)--(e\l)--cycle;}
\foreach \k in {0,2} {\draw[-{Stealth[length=2mm]},thick] (0,0)--(e\k);}
\foreach \k in {1,3,4} {\draw[{Stealth[length=2mm]}-,thick] (0,0)--(e\k);}
\foreach \k in {0,...,4} {\node[font=\scriptsize] at ({126+72*\k}:1.5) {$\k$};}
\fill (0,0) circle (1.6pt);
\end{scope}
\begin{scope}[xshift=11.2cm]
\foreach \k in {0,...,4} {\coordinate (n\k) at ({126+72*\k}:0.9);}
\draw[gray!60] (n0)--(n1)--(n2)--(n3)--(n4)--cycle;
\foreach \k in {0,2} {\node[vx] at (n\k) {};}
\foreach \k in {1,3,4} {\node[ov] at (n\k) {};}
\foreach \k in {0,...,4} {\node[font=\scriptsize] at ({126+72*\k}:1.2) {$\k$};}
\end{scope}
\node[font=\small] at (9.4,-2.2) {(b) $D(v)=\{0,2\}$};
\end{tikzpicture}
\caption{Descending links in a two-dimensional example, in which $W$ is the reflection group of a right-angled pentagon with sides $0,\dots,4$ and the nerve $L$ is a $5$-cycle. Each panel shows the star of a vertex $w$ of the Davis complex, with arrows pointing towards the lower endpoint of each edge, and the link of $w$, in which the directions in $D(v)$ are drawn black. A square lies in the descending link exactly when both of its directions are descending (dark square). (a) If $D(v)=\{0,1\}$, the descending link is an edge. (b) If $D(v)=\{0,2\}$, the descending link consists of two points.}\label{fig:links}
\end{figure}

Figure~\ref{fig:links} illustrates descending links in dimension $2$.

\subsection{Blocks and incoming boundaries}\label{sec:blocks}
Let $\Pi\subset I$ be a set of pairwise disjoint facets (\emph{ports}), and put
\[
W^\Pi=\langle r_i:\ i\in I\setminus\Pi\rangle,\qquad V^\Pi=\rho(W^\Pi),\qquad\Gamma^\Pi=\Gamma\cap W^\Pi .
\]
The hyperplanes containing $wF_i$ ($w\in W$, $i\in\Pi$) cut $\Hy^7$ into convex regions, and the closure $\widetilde N$ of the region containing $\operatorname{int}P$ is $\bigcup_{w\in W^\Pi}wP$. The \emph{block} of $(P,\Pi)$ is
\[
N=N(P,\Pi)=\widetilde N/\Gamma^\Pi .
\]
It is a finite-volume hyperbolic manifold with totally geodesic boundary, tessellated by $|V^\Pi|$ copies of $P$, with $\pi_1(N)=\Gamma^\Pi$ (and $N=\Hy^7/\Gamma$ when $\Pi=\emptyset$). Its character is $\varphi^\Pi=\varphi|_{\Gamma^\Pi}$, which depends only on the states and weights of $I\setminus\Pi$. A boundary facet $wF_i$ ($w\in W^\Pi$, $i\in\Pi$) is \emph{incoming} if $i\in D(\rho(w))$. By \eqref{eq:colors}, crossing facets that meet $F_i$ does not change whether $i\in D(v)$, so each boundary hyperplane of $\widetilde N$ is either wholly incoming or wholly outgoing. Let $\partial^-\widetilde N$ be the union of the incoming boundary hyperplanes and $\partial^-N=\partial^-\widetilde N/\Gamma^\Pi$ the \emph{incoming boundary}.

The Davis complex $\Sigma^\Pi\subset\Sigma$ of $W^\Pi$ consists of all cubes $wW_S$ with $w\in W^\Pi$ and $S\cap\Pi=\emptyset$. Its descending and ascending links at $w$ ($v=\rho(w)$) are the full subcomplexes of $L$ spanned by $D(v)\setminus\Pi$ and $D(v+\one)\setminus\Pi$. Similarly, for $i\in I$ let $I_{(i)}=\{j\in I:F_j\cap F_i\ne\emptyset\}$ and $W_{(i)}=\langle r_j:j\in I_{(i)}\rangle$. For any $u\in W$, the character
\[
\varphi_{(i)}^u\colon\Gamma\cap W_{(i)}\to\Z,\qquad\varphi_{(i)}^u(\gamma)=\varphi(u\gamma u^{-1}),
\]
admits the equivariant Morse function $w\mapsto f(uw)$ on the Davis complex of $W_{(i)}$, with descending and ascending links spanned by $D(v)\cap I_{(i)}$ and $D(v+\one)\cap I_{(i)}$ ($v=\rho(uw)$). Applying \cite[Theorem~4.1(3)]{BB97} (rescaling $\varphi$ and $f$ by the index of $\operatorname{im}\varphi$ in $\Z$ when nonzero, while the zero-character case is immediate since the quotient Davis complex is finite) yields:

\begin{lemma}\label{lem:blockfp}
If the full subcomplex of $L$ spanned by $D(v)\setminus\Pi$ (resp.\ $D(v)\cap I_{(i)}$) is simply connected for every $v\in V$, then $\ker\varphi^\Pi$ (resp.\ $\ker\varphi_{(i)}^u$ for all $u\in W$) is finitely presented.
\end{lemma}

\subsection{Explicit states, weights, and the five blocks}\label{sec:data}
On $P^7$ we choose the state\footnote{Found by a computer search minimizing $1000\sum_v(b_0(L_v)-1)+100\sum_v b_1(L_v)+\sum_v b_2(L_v)$ over $\F_2$. Every link of this state is simply connected (Proposition~\ref{prop:links}).} $\state_P\in\F_2^{56}$ given by
\[
\state_P=\texttt{10111101 01001011 10011011 00010000 11011101 10111011 01010011},
\]
reading bit $r$ of block $n$ ($0\le r\le7$, $0\le n\le6$) as the coordinate of $F_{8n+r}$. The weights $a^P_j$ for $j\ne8$ are listed in Table~\ref{tab:capweights}, while the weight $a^P_8>0$, needed only when $m=0$, will be fixed in Section~\ref{sec:novikov}.

On $Q$ we take the state\footnote{Found by a computer search requiring agreement with $\state_P$ on facets meeting $E_0$ (Lemma~\ref{lem:compatible}) and opposite values on $E_i$ and $\jmath(E_i)$ for every $i$.} $\state_Q\in\F_2^{182}$ recorded in Appendix~\ref{app:data}, together with the fixed-facet weights $a_i$ ($i\in\If$) in Table~\ref{tab:Qweights}. For the paired facets we set
\[
a_i=L\quad(i\in\Ip\setminus\{0,181\}),\qquad a_{181}=a_*,
\]
where the positive integers $L$ and $a_*$ will be chosen in Section~\ref{sec:novikov}. The weight of $E_0$ is never needed. Direct inspection of Appendix~\ref{app:data} verifies:

\begin{lemma}\label{lem:compatible}
Whenever $j\sim8$ and $g_i=e_j$ (Lemma~\ref{lem:Q}(5)), $\state_P$ and $\state_Q$ agree at $F_j$ and $E_i$, and $a^P_j=a_i$. Moreover, $(\state_P)_8=0$, $(\state_Q)_0=1$, and $(\state_Q)_{181}=0$.
\end{lemma}

On $S_0=P^7$ we use $(\state_P,a^P)$. On $S_m$ ($m\ge1$), each facet $F_j$ of $P^7$ ($j\ne8$) retains its state and weight $a^P_j$, while each facet $h_k(E_i)$ of $Q_k$ receives weight $a_i$ and state $(\state_Q)_i$ ($i\in\If$) or $(\state_Q)_i+k-1\pmod2$ ($i\in\Ip$). Thus $Q_k$ carries $\state_Q$ for $k$ odd and $\state_Q^\sharp$ (obtained by flipping all paired-facet coordinates) for $k$ even. By Lemma~\ref{lem:compatible}, these assignments agree on merged facets of consecutive pieces, defining a character $\alpha\colon\Gamma_m\to\Z$.

We work with five blocks:
\begin{center}
\begin{tabular}{lllll}
\toprule
block&polytope&ports $\Pi$&group&character\\\midrule
$M_0$&$P^7$&$\emptyset$&$\Gamma_0$&$\alpha$\\
$C$ (cap)&$P^7$&$\{8\}$&$\Gamma_C$&$\varphi_C$\\
$B$ (bridge)&$Q$&$\{0,181\}$&$\Gamma_B$&$\varphi_B$\\
$E$ (end block)&$Q$&$\{0\}$&$\Gamma_E$&$\varphi_E$\\
$X$ (core)&$Q$&$\Ip$&$\Gamma_X$&$\varphi_X$\\\bottomrule
\end{tabular}
\end{center}
with data $(\state_P,a^P)$ for $M_0,C$ and $(\state_Q,a)$ for $B,E,X$. The non-port colors span $V$ for $M_0,C,B,E$ (tessellated by $256$ tiles) and $U$ for $X$ (tessellated by $128$ tiles, see Appendix~\ref{app:computations}).

\subsection{Link topology and cusp periods}

\begin{proposition}\label{prop:links}
\begin{enumerate}
\item For $(P^7,\state_P)$ and every $v\in V$, the complex $L_v$ is simply connected, and
\[
\widetilde H_k(L_v;\Z)\cong
\begin{cases}
\Z,& k=2\text{ and }v\in\{y_1,\dots,y_{24}\},\\
\Z,& k=3\text{ and }v\in\{x_1,\dots,x_{24}\},\\
0,& \text{otherwise,}
\end{cases}
\]
with $x_j,y_i$ listed in Table~\ref{tab:pivots}. The full subcomplexes spanned by $D(v)\setminus\{8\}$ and $D(v)\cap I_{(8)}$ are also simply connected.
\item For $(Q,\state_Q)$ and every $v\in V$, the complex $L_v$ is simply connected. Its reduced homology groups $\widetilde H_2(L_v;\Z)$ and $\widetilde H_3(L_v;\Z)$ are free abelian of the same rank $k(v)\le4$, and $\widetilde H_k(L_v;\Z)=0$ for $k\ne2,3$. The number of $v\in V$ with $k(v)=0,1,2,3,4$ is $38,88,68,44,18$, so $\sum_{v\in V}k(v)=428$ and $\sum_{v\in U}k(v)=214$. The full subcomplexes spanned by $D(v)\setminus\{0\}$, $D(v)\setminus\{0,181\}$, $D(v)\cap I_{(0)}$, and $D(v)\cap I_{(181)}$ are simply connected as well.
\end{enumerate}
\end{proposition}

\begin{proof}
Let $K$ be one of these subcomplexes, and pick a spanning tree $T$ of $K^{(1)}$. The edge-path group of $K$ has the presentation
\[
\pi_1(K)\cong\bigl\langle\,g_{ij}\ (\{i<j\}\in K)\ \bigm|\ g_{ij}=1\ (\{i,j\}\in T),\ \ g_{ij}\,g_{jk}=g_{ik}\ (\{i<j<k\}\in K)\,\bigr\rangle.
\]
Explicit Tietze transformations reduce this presentation to the trivial one, so $\pi_1(K)=1$. For the homology of $L_v$, we repeatedly split off acyclic summands $\Z\sigma\oplus\Z\partial\sigma$ with $\langle\partial\sigma,\tau\rangle=\pm1$ as in Section~\ref{sec:local}. What remains is a chain-homotopy-equivalent complex with zero differential, concentrated in the stated degrees and ranks. Both the Tietze reductions and the cancellation pairs are checked by computer (Appendix~\ref{app:computations}).
\end{proof}

\begin{corollary}\label{cor:fp}
The kernels of $\varphi_C,\varphi_B,\varphi_E$, of $\alpha$ on $\Gamma_0$, and of $\varphi_{(8)}^u$ on $P^7$ and $\varphi_{(0)}^u,\varphi_{(181)}^u$ on $Q$ (for all $u$) are finitely presented.
\end{corollary}

\begin{proof}
This follows from Lemma~\ref{lem:blockfp} and Proposition~\ref{prop:links}.
\end{proof}

Let $(F_\alpha,F_\beta)$ be a tangent pair of $P$, and let $w\in W$ with $v=\rho(w)$. Then $c_\alpha\cdot c_\beta=1$, so $(r_\alpha r_\beta)^2\in\Gamma$. Summing the height increments along $w,wr_\alpha,wr_\alpha r_\beta,wr_\alpha r_\beta r_\alpha,w(r_\alpha r_\beta)^2$ gives
\begin{equation}\label{eq:period}
\varphi\bigl(w(r_\alpha r_\beta)^2w^{-1}\bigr)=2\bigl(\alpha(v,\alpha)-\alpha(v,\beta)\bigr).
\end{equation}
This is nonzero whenever $a_\alpha\ne a_\beta$ or $|\{\alpha,\beta\}\cap D(v)|=1$. A finite computation (Appendix~\ref{app:computations}) gives the next lemma.

\begin{lemma}\label{lem:peripheral}
\begin{enumerate}
\item For every ideal vertex $\xi$ of $Q$ and every $v\in V$, at least one of the five tangent fixed pairs $(E_\alpha,E_\beta)$ at $\xi$ satisfies $\alpha(v,\alpha)\ne\alpha(v,\beta)$ for $\state_Q$ and Table~\ref{tab:Qweights}.
\item For every ideal vertex $\xi$ of $P^7$ and every $v\in V$, at least one tangent pair $(F_\alpha,F_\beta)$ at $\xi$ not containing $F_8$ satisfies $\alpha(v,\alpha)\ne\alpha(v,\beta)$ for $\state_P$ and Table~\ref{tab:capweights}.
\item The character $\varphi_C\colon\Gamma_C\to\Z$ is surjective, because the fundamental cycle periods of a spanning tree of the Schreier graph on $V$ have greatest common divisor $1$.
\end{enumerate}
\end{lemma}

\section{Homological perturbation and relative Novikov acyclicity}\label{sec:novikov}

\subsection{Novikov rings, weighted norms, and perturbation}
Given a group $G$ and a homomorphism $\chi\colon G\to\Z$, the \emph{Novikov ring} $\Nov{G}{\chi}$ is the ring of formal sums $x=\sum_{g\in G}x_g\,g$ with $x_g\in\Z$ whose support intersects $\{g:\chi(g)>c\}$ in a finite set for every $c\in\R$. Writing $\nu_\chi(x)=\max\{\chi(g):x_g\ne0\}$ (and $\nu_\chi(0)=-\infty$), we have $\nu_\chi(xy)\le\nu_\chi(x)+\nu_\chi(y)$ and $\nu_\chi(x+y)\le\max\{\nu_\chi(x),\nu_\chi(y)\}$. Any subgroup $G'<G$ gives a subring $\Nov{G'}{\chi}\subseteq\Nov{G}{\chi}$.

Set $R=\Nov{G}{\chi}$. By a \emph{module with heights} we mean a finitely generated free left $R$-module with a chosen basis $(\mathbf b_k)$ and an integer height $h(\mathbf b_k)$ for each basis element. Given an $R$-linear map $A(\mathbf b_k)=\sum_i A_{ik}\mathbf b'_i$ between two such modules, we set
\[
\|A\|=\max_{i,k}\bigl(\nu_\chi(A_{ik})+h(\mathbf b'_i)-h(\mathbf b_k)\bigr).
\]
Saying $\|A\|\le-c$ means that every term of $A$ drops the height by at least $c$. Clearly $\|A\circ A'\|\le\|A\|+\|A'\|$ and $\|A+A'\|\le\max\{\|A\|,\|A'\|\}$.

\begin{lemma}\label{lem:neumann}
If $E$ is an endomorphism of a module with heights and $\|E\|<0$, then $1+E$ is invertible with $\|(1+E)^{-1}\|\le0$.
\end{lemma}

\begin{proof}
Since $\|E^n\|\le-n$ and only finitely many basis heights appear, the Neumann series $\sum_{n\ge0}(-E)^n$ converges entrywise in $R$ to $(1+E)^{-1}$, with norm at most $0$.
\end{proof}

A bounded complex of finitely generated free $R$-modules is acyclic if and only if it is contractible, meaning there is a degree-$1$ map $\Theta$ with $d\Theta+\Theta d=1$. Such a complex stays acyclic under base change along any ring map $R\to R'$.

\begin{lemma}\label{lem:perturb}
Suppose $(C,d)$ is a bounded complex of modules with heights, $\Theta$ is a contraction of $(C,d)$, and $\delta\colon C\to C$ is a degree-$(-1)$ map with $(d+\delta)^2=0$. If $\|\delta\|+\|\Theta\|<0$, then $(C,d+\delta)$ is acyclic.
\end{lemma}

\begin{proof}
Put $d'=d+\delta$ and $u=d'\Theta+\Theta d'=1+\delta\Theta+\Theta\delta$. Since $\|\delta\Theta+\Theta\delta\|<0$, Lemma~\ref{lem:neumann} implies that $u$ is a degreewise isomorphism. Because $d'^2=0$, $u$ is a chain automorphism of $(C,d')$ that is chain homotopic to $0$. Therefore $H_*(C,d')=0$.
\end{proof}

\subsection{Cellular models of blocks and general invariance principles}\label{sec:blockcomplex}
For $(P,\Pi)$ as in Section~\ref{sec:blocks}, let $\Sigma^+\subset\Sigma$ be the subcomplex formed by the cubes $wW_S$ with $w\in W^\Pi$ and $S\cap\Pi\subseteq D(\rho(w))$, and let $\Sigma^-\subset\Sigma^+$ be the subcomplex of cubes with no vertex in $W^\Pi$. Since no clique contains two ports, $\Sigma^+$ is $\Sigma^\Pi$ with a product collar attached along each incoming boundary hyperplane, and $\Sigma^-$ is the union of the outer ends of these collars. We set $C(P,\Pi)=C_*(\Sigma^+,\Sigma^-)$. This is a bounded complex of finitely generated free left $\Z\Gamma^\Pi$-modules, reducing to $C_*(\Sigma)$ when $\Pi=\emptyset$.

\begin{lemma}\label{lem:model}
For any $\Gamma^\Pi$-invariant CW structure on $\widetilde N$ having $\partial^-\widetilde N$ as a subcomplex, $C(P,\Pi)$ is chain homotopy equivalent to $C_*(\widetilde N,\partial^-\widetilde N)$ over $\Z\Gamma^\Pi$.
\end{lemma}

\begin{proof}
Consider the two free $\Gamma^\Pi$-CW pairs $(\Sigma^+,\Sigma^-)$ and $(\widetilde N,\partial^-\widetilde N)$. Here $\Sigma^+$ deformation retracts onto the contractible complex $\Sigma^\Pi$, while $\widetilde N$ is convex. The components of $\Sigma^-$ are the translated Davis complexes $wr_i\Sigma_{(i)}$ of the facet subgroups $W_{(i)}$, which are contractible and correspond $\Gamma^\Pi$-equivariantly to the convex boundary components $wF_i\cap\widetilde N$ of $\partial^-\widetilde N$. By equivariant cellular induction, the two pairs are $\Gamma^\Pi$-equivariantly homotopy equivalent, so their cellular chain complexes are chain homotopy equivalent over $\Z\Gamma^\Pi$.
\end{proof}

\begin{definition}\label{def:acyclic}
We call the block of $(P,\Pi)$ \emph{acyclic} if $H_*(\Nov{\Gamma^\Pi}{\varphi^\Pi}\otimes_{\Z\Gamma^\Pi}C(P,\Pi))=0$, or equivalently $H_*(N,\partial^-N;\Nov{\Gamma^\Pi}{\varphi^\Pi})=0$.
\end{definition}

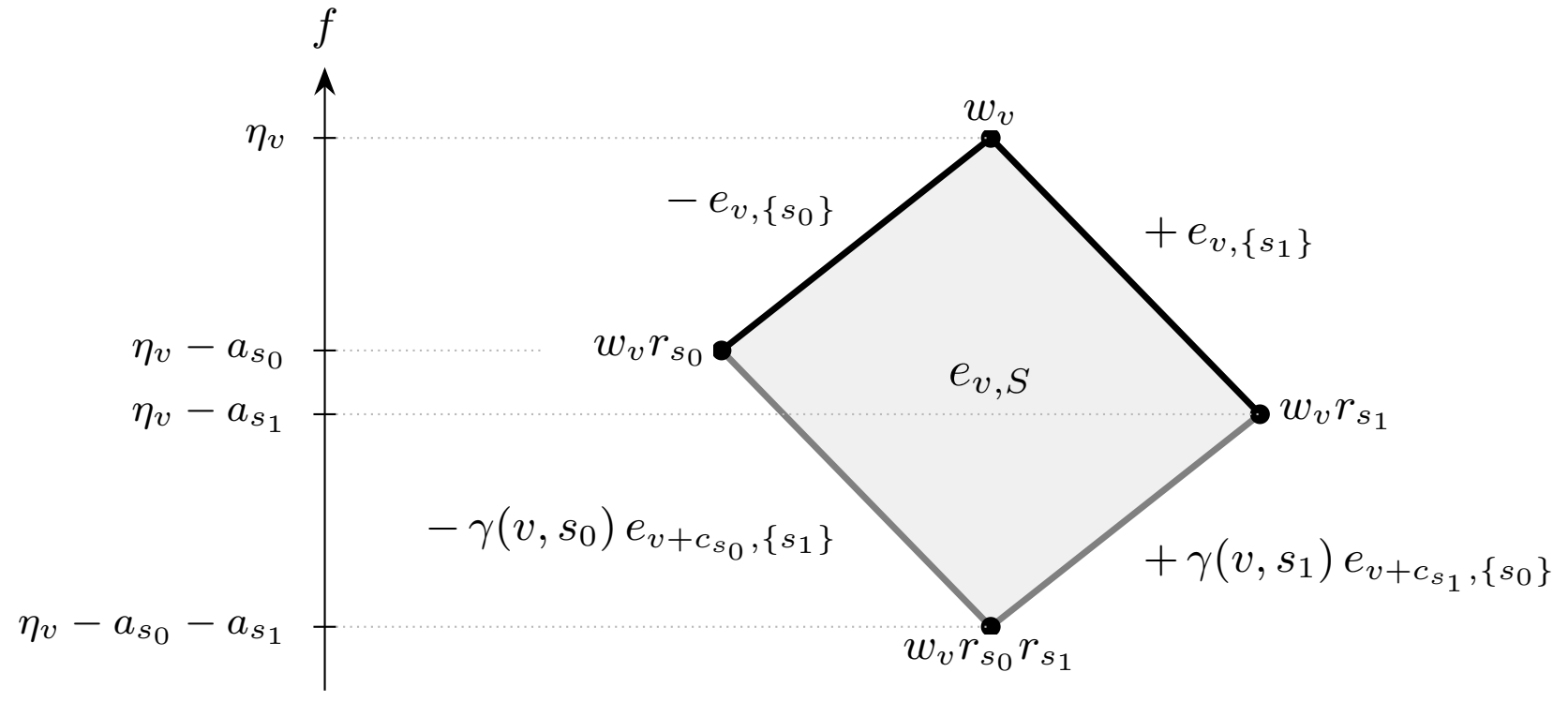
\begin{figure}[ht]
\centering
\begin{tikzpicture}[x=1cm,y=1cm,vx/.style={circle,fill=black,inner sep=1.4pt},lab/.style={font=\footnotesize}]
\coordinate (T) at (0,0); \coordinate (L) at (-1.9,-1.5); \coordinate (R) at (1.9,-1.95); \coordinate (Bt) at (0,-3.45);
\fill[gray!12] (T)--(L)--(Bt)--(R)--cycle;
\node[font=\small] at (0,-1.72) {$e_{v,S}$};
\draw[line width=1.3pt] (T)--(R) node[lab,midway,above right=-1pt] {$+\,e_{v,\{s_1\}}$};
\draw[line width=1.3pt] (T)--(L) node[lab,midway,above left=-1pt] {$-\,e_{v,\{s_0\}}$};
\draw[line width=1.3pt,gray] (L)--(Bt) node[lab,midway,below left=-1pt,text=black] {$-\,\gamma(v,s_0)\,e_{v+c_{s_0},\{s_1\}}$};
\draw[line width=1.3pt,gray] (R)--(Bt) node[lab,midway,below right=-1pt,text=black] {$+\,\gamma(v,s_1)\,e_{v+c_{s_1},\{s_0\}}$};
\foreach \P in {T,L,R,Bt} \node[vx] at (\P) {};
\node[lab,above,fill=white,inner sep=1pt] at ($(T)+(0,0.05)$) {$w_v$};
\node[lab,left,fill=white,inner sep=1pt] at ($(L)+(-0.05,0)$) {$w_vr_{s_0}$};
\node[lab,right] at (R) {$w_vr_{s_1}$};
\node[lab,below,fill=white,inner sep=1pt] at ($(Bt)+(0,-0.05)$) {$w_vr_{s_0}r_{s_1}$};
\draw[-{Stealth[length=2mm]}] (-4.7,-3.9)--(-4.7,0.5) node[above,font=\footnotesize] {$f$};
\foreach \y/\t in {0/{\eta_v},-1.5/{\eta_v-a_{s_0}},-1.95/{\eta_v-a_{s_1}},-3.45/{\eta_v-a_{s_0}-a_{s_1}}}
  {\draw (-4.78,\y)--(-4.62,\y); \node[font=\scriptsize,anchor=east] at (-4.82,\y) {$\t$};}
\draw[gray!55,densely dotted] (-4.62,0)--(T) (-4.62,-1.95)--(R) (-4.62,-3.45)--(Bt) (-4.62,-1.5)--(-3.15,-1.5);
\end{tikzpicture}
\caption{A square $e_{v,S}$ with $S=\{s_0<s_1\}$, drawn with the height increasing upwards. The two upper sides contain the top vertex $w_v$ and form $d_0\,e_{v,S}$ (black). The two lower sides form $d_1\,e_{v,S}$ (gray). If $s_0$ is a port, the lower side in the direction $s_0$ lies in $\Sigma^-$, and the corresponding term is omitted.}\label{fig:faces}
\end{figure}

Each cube of $\Sigma^+\setminus\Sigma^-$ is uniquely of the form $wW_S$ with $w\in W^\Pi$ and $S\subseteq D(\rho(w))$, where $w$ is the unique vertex at which $f$ attains its maximum on the cube (its \emph{top} vertex). Pick a vertex $w_v\in W^\Pi$ with $\rho(w_v)=v$ for every $v\in V^\Pi$, and write $\eta_v=f(w_v)$. Then the cubes $e_{v,S}=w_vW_S$, with $v\in V^\Pi$ and $S\subseteq D(v)$ a clique, give a $\Z\Gamma^\Pi$-basis of $C(P,\Pi)$ in degree $|S|$. Ordering $S=\{s_0<\dots<s_{n-1}\}$, the cellular boundary splits as $\partial=d_0+d_1$ (Figure~\ref{fig:faces}), with
\begin{align*}
d_0\,e_{v,S}&=\sum_{k=0}^{n-1}(-1)^k\,e_{v,S\setminus s_k},\\
d_1\,e_{v,S}&=\sum_{k:\ s_k\notin\Pi}(-1)^{k+1}\,\gamma(v,s_k)\,e_{v+c_{s_k},\,S\setminus s_k},\qquad \gamma(v,s)=w_v r_s w_{v+c_s}^{-1}\in\Gamma^\Pi.
\end{align*}
Deleting from $d_1$ the lower faces in any subset of facet directions leaves a map whose sum with $d_0$ still squares to $0$, as canceling pairs in $\partial^2 e_{v,S}$ use lower faces in identical directions. Assigning $e_{v,S}$ height $\eta_v$ makes $\Nov{\Gamma^\Pi}{\varphi^\Pi}\otimes C(P,\Pi)$ a complex of modules with heights: $d_0$ preserves height, while the $s$-term of $d_1$ lowers height by $a_s$.

Given any state $\state$, put $\state^\sharp_i=\state_i+c_i\cdot\eps$, which agrees with Section~\ref{sec:data} on $Q$.

\begin{lemma}\label{lem:sharp}
If $\eps\in V^\Pi$, then $(P,\Pi)$ is acyclic for $\state^\sharp$ if and only if it is acyclic for $\state$ with the same weights. Moreover, the kernels of the two characters are isomorphic.
\end{lemma}

\begin{proof}
Pick $u\in W^\Pi$ with $\rho(u)=\eps$. Because $D^\sharp(v)=D(v+\eps)$, the height functions satisfy $f^\sharp(w)=f(uw)-f(u)$ on $W^\Pi$, so $\varphi^\sharp(\gamma)=\varphi(u\gamma u^{-1})$. Left translation by $u$ carries $(\Sigma^{+\sharp},\Sigma^{-\sharp})$ isomorphically onto $(\Sigma^+,\Sigma^-)$, intertwining the $\Gamma^\Pi$-actions via conjugation by $u$. This conjugation extends to an isomorphism between the two Novikov rings and identifies the two complexes.
\end{proof}

\begin{lemma}\label{lem:restore}
Let $\Pi'\subset\Pi$ be sets of ports of $P$ such that $V^{\Pi'}=V^\Pi$. If $(P,\Pi)$ is acyclic, then there exists $M_0$ such that $(P,\Pi')$ is acyclic whenever $a_s\ge M_0$ for all $s\in\Pi\setminus\Pi'$.
\end{lemma}

\begin{proof}
Choose representatives $w_v\in W^\Pi$ for all $v\in V^\Pi=V^{\Pi'}$. As a graded module, $C(P,\Pi')=\Z\Gamma^{\Pi'}\otimes_{\Z\Gamma^\Pi}C(P,\Pi)$, and its differential is $d+\delta$, where $d$ is induced from $C(P,\Pi)$ and $\delta$ consists of the lower faces in directions $s\in\Pi\setminus\Pi'$. Neither the heights $\eta_v=f(w_v)$ nor the restriction $\varphi^{\Pi'}|_{\Gamma^\Pi}=\varphi^\Pi$ depends on the weights $\{a_s:s\in\Pi\setminus\Pi'\}$. A contraction $\Theta$ of $\Nov{\Gamma^\Pi}{\varphi^\Pi}\otimes C(P,\Pi)$ induces a contraction of $\Nov{\Gamma^{\Pi'}}{\varphi^{\Pi'}}\otimes C(P,\Pi)$ with the same norm $\|\Theta\|$, which is independent of these weights. Meanwhile $\|\delta\|\le-\min_{s\in\Pi\setminus\Pi'}a_s$, so Lemma~\ref{lem:perturb} gives acyclicity whenever every $a_s>\|\Theta\|$.
\end{proof}

\subsection{Local strong deformation retractions and the dominance criterion}\label{sec:local}
Setting $C(v)=\bigoplus_{S\subseteq D(v)}\Z\,e_{v,S}$, we have
\begin{equation}\label{eq:Csplit}
\bigl(C(P,\Pi),d_0\bigr)=\bigoplus_{v\in V^\Pi}\Z\Gamma^\Pi\otimes_{\Z}\bigl(C(v),d_0\bigr),
\end{equation}
and $e_{v,S}\mapsto S$ identifies $(C(v),d_0)$ with the shifted augmented simplicial chain complex $\widetilde C_{*-1}(L_v)$, with homology $H_k(v)\cong\widetilde H_{k-1}(L_v;\Z)$.

We contract each $(C(v),d_0)$ onto $(H(v),0)$ via a \emph{strong deformation retraction} $(\iota_v,\pi_v,\kappa_v)$ satisfying
\begin{equation}\label{eq:sdr}
\begin{gathered}
d_0\,\iota_v=0,\qquad \pi_v\,d_0=0,\qquad \pi_v\,\iota_v=1,\qquad d_0\,\kappa_v+\kappa_v\,d_0=1-\iota_v\,\pi_v,\\
\kappa_v^2=0,\qquad \kappa_v\,\iota_v=0,\qquad \pi_v\,\kappa_v=0.
\end{gathered}
\end{equation}
Suppose $\sigma$ and $\tau$ are basis elements with $\langle d_0\sigma,\tau\rangle=\theta=\pm1$. Let $C'(v)$ be the span of the remaining basis elements and define
\[
\iota'(x)=x-\theta\langle d_0x,\tau\rangle\sigma,\qquad
\pi'(x)=x-\langle x,\sigma\rangle\sigma-\theta\langle x,\tau\rangle d_0\sigma,\qquad
\kappa'(x)=\theta\langle x,\tau\rangle\sigma,
\]
with $d_0'=\pi'd_0\iota'$. This step splits off the acyclic summand $\Z\sigma\oplus\Z d_0\sigma$ (\emph{cancelling} $(\sigma,\tau)$). For our polytopes and states, $N$ iterations reduce the differential to $0$ (Proposition~\ref{prop:links}). Composing these steps by
\[
\iota_v=\iota^{(1)}\cdots\iota^{(N)},\qquad
\pi_v=\pi^{(N)}\cdots\pi^{(1)},\qquad
\kappa_v=\sum_{n=1}^{N}\iota^{(1)}\cdots\iota^{(n-1)}\,\kappa^{(n)}\,\pi^{(n-1)}\cdots\pi^{(1)}
\]
gives a strong deformation retraction satisfying \eqref{eq:sdr} \cite{Skoldberg06}, determined by the cancellation order in Appendix~\ref{app:computations}. Taking the direct sum over $v\in V^\Pi$ gives a height-preserving strong deformation retraction $(\iota,\pi,\kappa)$ from $(C(P,\Pi),d_0)$ onto $(H,0)$, where $H=\bigoplus_v\Z\Gamma^\Pi\otimes H(v)$.

Write $R=\Nov{\Gamma^\Pi}{\varphi^\Pi}$, $\widehat C=R\otimes C(P,\Pi)$, and $\widehat H=R\otimes H$. Because $\|(\kappa d_1)^n\|\le-n$, the operator $1+\kappa d_1$ is invertible over $R$ by Lemma~\ref{lem:neumann}. Applying the basic perturbation lemma \cite{Brown65,Crainic04}, we obtain a chain homotopy equivalence between $(\widehat C,d_0+d_1)$ and $(\widehat H,D)$, where
\[
D=\pi\,d_1\,(1+\kappa d_1)^{-1}\iota=\sum_{r\ge1}(-1)^{r-1}D^{(r)},\qquad D^{(r)}=\pi\,d_1(\kappa d_1)^{r-1}\iota .
\]
Suppose now that $H$ lives only in degrees $3$ and $4$ and that $H_3$ and $H_4$ are free of the same rank $n$, which is the case for $(P^7,\state_P)$ and $(Q,\state_Q)$ by Proposition~\ref{prop:links}. Pick bases $\mathbf f_1,\dots,\mathbf f_n$ of $H_3$ with $\mathbf f_i\in H_3(y_i)$ and $\mathbf e_1,\dots,\mathbf e_n$ of $H_4$ with $\mathbf e_j\in H_4(x_j)$, endowed with heights $\eta_{y_i}$ and $\eta_{x_j}$. The differential is then a single $n\times n$ matrix $D\colon\widehat H_4\to\widehat H_3$ over $R$, and
\begin{equation}\label{eq:criterion}
\text{the block is acyclic if and only if $D$ is an isomorphism.}
\end{equation}

Call a sequence $p=(F_1,\dots,F_r)$ of facets in $I\setminus\Pi$ a \emph{path} from $x\in V^\Pi$. Setting $v_0=x$ and $v_k=v_{k-1}+c_{F_k}$, the path terminates at $v_r$, and we attach to it the group element and signed weight
\[
\gamma_p=w_x\,r_{F_1}\cdots r_{F_r}\,w_{v_r}^{-1}\in\Gamma^\Pi,\qquad a(p)=-\sum_{k=1}^r\alpha(v_{k-1},F_k).
\]
Lemma~\ref{lem:cocycle} gives
\begin{equation}\label{eq:phipath}
\varphi^\Pi(\gamma_p)=\eta_x-\eta_{v_r}-a(p),
\end{equation}
so $a(p)$ is determined by $x$, $v_r$, and $\gamma_p$, and adds when paths are concatenated. When $F_k\in D(v_{k-1})\setminus\Pi$ for all $k$, we call $w=p$ a \emph{descending word}, and then $a(w)=\sum_{k=1}^r a_{F_k}$. No descending word can be shortened in $W$: by the deletion condition, a shortening would match two copies of some generator $r_F$ with only commuting letters between them, forcing the second crossing of $F$ to go upwards. Consequently, two descending words represent the same element of $W$ if and only if they differ by commuting adjacent letters. Since $\iota,\pi,\kappa$ do not change vertex labels while $d_1$ steps along descending facets in $I\setminus\Pi$, each entry of $D^{(r)}$ has the form $D^{(r)}_{ij}=\sum_w n_w\gamma_w$ ($n_w\in\Z$), summed over descending words $w$ of length $r$ from $x_j$ to $y_i$.

To prove that $D$ is an isomorphism, we allow a unitriangular change of basis on $\widehat H_4$ and shift the heights of the target basis elements $\mathbf f_i$ by integers $u_i$, so that a monomial matrix $Q_0$ becomes the strictly dominant part of $D$ while all remaining terms of orders $1$ and $2$ and all higher-order terms ($r\ge3$) strictly lower the shifted height.

\begin{proposition}\label{lem:dominance}
Suppose $(\mathbf e'_j)_{j=1}^n$ is a basis of $\widehat H_4$ of the form $\mathbf e'_j=\mathbf e_j+\sum\pm\gamma_p\mathbf e_k$, where each $p$ is a path from $x_j$ to $x_k$ ($k\ne j$) and the relation $k\prec j$ ($\mathbf e_k$ appears in $\mathbf e'_j-\mathbf e_j$) has no cycles. Write $D(\mathbf e'_j)=\sum_i D'_{ij}\mathbf f_i$, and let $D'^{\,\le2}_{ij}\in\Z\Gamma^\Pi$ be obtained by replacing $D$ with $D^{(1)}-D^{(2)}$ in
\begin{equation}\label{eq:Dprime}
D'_{ij}=D_{ij}+\sum\pm\gamma_p\,D_{ik}.
\end{equation}
Suppose there exist a bijection $\beta$ of $\{1,\dots,n\}$, paths $w_j$ from $x_j$ to $y_{\beta(j)}$, and integers $u_1,\dots,u_n$ such that:
\begin{enumerate}
\item for every $j$, the coefficient of $\gamma_{w_j}$ in $D'^{\,\le2}_{\beta(j)j}$ is $\epsilon_j=\pm1$;
\item every other term $n\,\gamma_q$ of $D'^{\,\le2}_{ij}$ (after collecting equal group elements) satisfies $a(q)-a(w_j)\ge u_i-u_{\beta(j)}+1$;
\item for all $i,j$ and every path $p$ occurring in $\mathbf e'_j-\mathbf e_j$ (including the empty path from $x_j$), $a(w_j)-a(p)+u_i-u_{\beta(j)}<3\,a_{\min}$, where $a_{\min}=\min_{F\notin\Pi}a_F$.
\end{enumerate}
Then $D$ is an isomorphism and the block is acyclic.
\end{proposition}

\begin{proof}
Because $\prec$ has no cycles, we can order the indices so that $(\mathbf e'_j)$ is obtained from $(\mathbf e_j)$ by a unitriangular matrix over $\Z\Gamma^\Pi$. Thus $(\mathbf e'_j)$ is a basis of $\widehat H_4$. Let $Q_0\colon\widehat H_4\to\widehat H_3$ be the isomorphism $Q_0(\mathbf e'_j)=\epsilon_j\gamma_{w_j}\mathbf f_{\beta(j)}$, and shift the height of each $\mathbf f_i$ to $\eta_{y_i}+u_i$. Writing $D-Q_0=(D'^{\,\le2}-Q_0)+D'^{\,\ge3}$ in the bases $(\mathbf e'_j)$ and $(\mathbf f_i)$, where $D'^{\,\ge3}$ sums the contributions from $D^{(r)}$ for $r\ge3$, consider the endomorphism $E=(D-Q_0)Q_0^{-1}$ of $\widehat H_3$.

Each term $n\,\gamma_q\mathbf f_i$ of $(D-Q_0)(\mathbf e'_j)$ gives a term $\pm n\,\gamma_{w_j}^{-1}\gamma_q$ of $E$ from $\mathbf f_{\beta(j)}$ to $\mathbf f_i$. By \eqref{eq:phipath}, the resulting change in shifted height is
\[
\varphi^\Pi(\gamma_{w_j}^{-1}\gamma_q)+(\eta_{y_i}+u_i)-(\eta_{y_{\beta(j)}}+u_{\beta(j)})=a(w_j)-a(q)+u_i-u_{\beta(j)}.
\]
For $D'^{\,\le2}-Q_0$ this quantity is at most $-1$ by~(1) and~(2). For $D'^{\,\ge3}$, every term has $q=pw$ where $p$ is as in~(3) and $w$ is a descending word of length $r\ge3$. Hence $a(q)=a(p)+a(w)\ge a(p)+3a_{\min}$, so the height change is negative by~(3). Therefore $\|E\|<0$. By Lemma~\ref{lem:neumann}, $D=(1+E)Q_0$ is an isomorphism, and the block is acyclic by \eqref{eq:criterion}.
\end{proof}

\begin{figure}[ht]
\centering
\begin{tikzpicture}[font=\small,x=1cm,y=1cm,>={Stealth[length=2mm]},vx/.style={circle,fill=black,inner sep=1.5pt},ov/.style={circle,draw=black,fill=white,inner sep=1.5pt}]
\begin{scope}
\draw[->] (0,-0.2)--(0,3.9) node[above] {$f$};
\coordinate (p0) at (1.0,3.3); \coordinate (p1) at (1.7,2.7); \coordinate (p2) at (2.4,3.05); \coordinate (p3) at (3.1,2.2);
\coordinate (p4) at (3.6,1.5); \coordinate (p5) at (4.1,0.85); \coordinate (p6) at (4.6,0.2);
\foreach \a/\b in {p0/p1,p1/p2,p2/p3} {\draw[thick] (\a)--(\b);}
\foreach \a/\b in {p3/p4,p4/p5,p5/p6} {\draw[thick,gray] (\a)--(\b);}
\draw[thick,->] (p0)--($(p0)!0.62!(p1)$); \draw[thick,->] (p1)--($(p1)!0.62!(p2)$); \draw[thick,->] (p2)--($(p2)!0.62!(p3)$);
\draw[thick,gray,->] (p3)--($(p3)!0.62!(p4)$); \draw[thick,gray,->] (p4)--($(p4)!0.62!(p5)$); \draw[thick,gray,->] (p5)--($(p5)!0.62!(p6)$);
\foreach \p in {p1,p2,p4,p5} {\fill (\p) circle (1.1pt);}
\node[vx] at (p0) {}; \node[vx] at (p3) {}; \node[ov] at (p6) {};
\node[above right=-1pt] at (p0) {$\mathbf e_j$};
\node[right=2pt] at (p3) {$\gamma_p\,\mathbf e_k$};
\node[right=2pt] at (p6) {$\gamma_{pw}\,\mathbf f_i$};
\node[font=\footnotesize] at (2.05,3.45) {$p$};
\node[font=\footnotesize] at (3.55,0.95) {$w$};
\foreach \y/\t in {3.3/{\eta_{x_j}},2.2/{\eta_{x_j}-a(p)},0.2/{\eta_{x_j}-a(p)-a(w)}}
  {\draw (-0.07,\y)--(0.07,\y); \node[font=\scriptsize,anchor=east] at (-0.1,\y) {$\t$};}
\draw[gray!60,densely dotted] (0.07,3.3)--(p0) (0.07,2.2)--(p3) (0.07,0.2)--(p6);
\end{scope}
\end{tikzpicture}
\caption{A term $n\,\gamma_{pw}\,\mathbf f_i$ of $D(\mathbf e'_j)$: a path $p$ from $x_j$ to $x_k$, which may go upwards, followed by a descending word $w$ from $x_k$ to $y_i$.}\label{fig:dominance}
\end{figure}
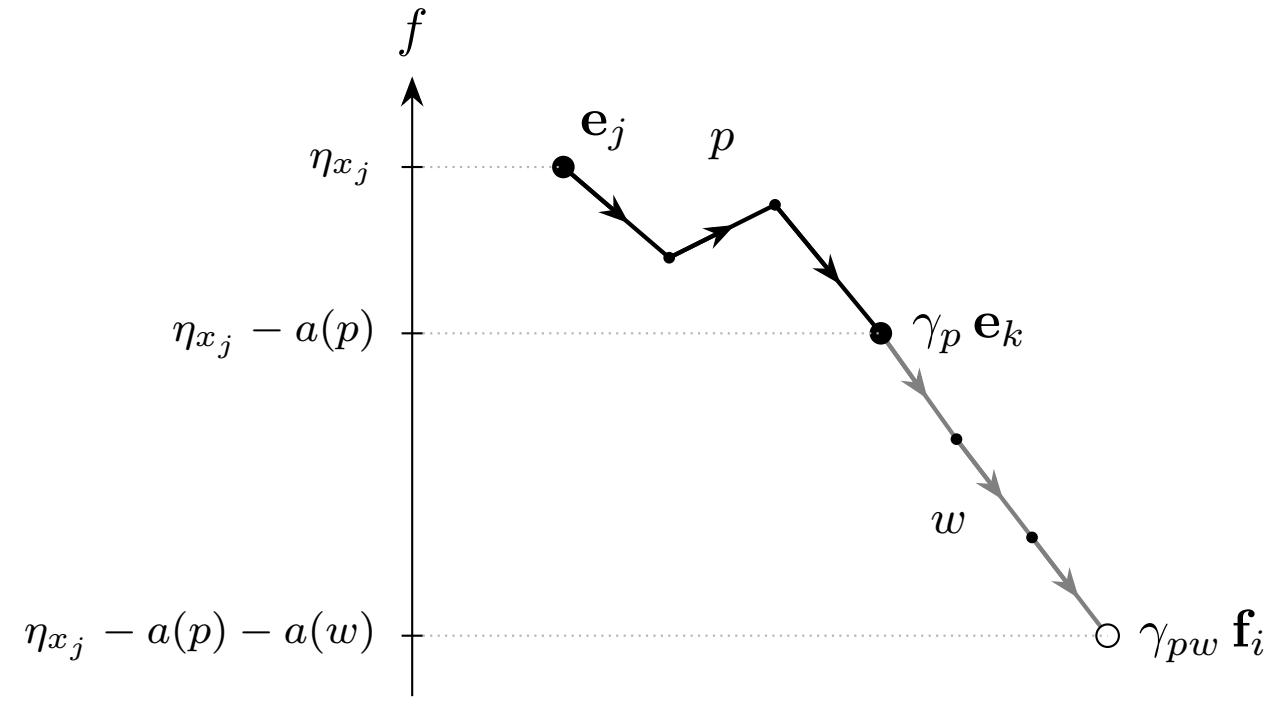

In Proposition~\ref{lem:dominance}, the paths $p$ and $w_j$ need not be descending. A step going upwards raises the height, and this is accounted for by the signed quantity $a(p)$ in condition~(3) (Figure~\ref{fig:dominance}).

\subsection{\texorpdfstring{Acyclicity of the cap $C$, $M_0$, and the core $X$}{Acyclicity of the cap C, M0, and the core X}}\label{sec:core}
For the cap $C=(P^7,\{8\})$ with state $\state_P$ and the weights of Table~\ref{tab:capweights}, we have $n=24$. Table~\ref{tab:pivots} lists $x_j,y_i$, a bijection $\beta$ of $\{1,\dots,24\}$, descending words $w_j$ of length $1$ or $2$ avoiding $r_8$, and integers $u_1,\dots,u_{24}$. A finite exact computation (Appendix~\ref{app:computations}, where $D^{(1)}$ has $16$ nonzero terms and $D^{(2)}$ has $163$) gives the following leading-term bounds.

\begin{proposition}\label{prop:capcert}
For the cap $C$, with $\mathbf e'_j=\mathbf e_j$ and the retractions of Appendix~\ref{app:computations}:
\begin{enumerate}
\item for every $j$, the coefficient of $\gamma_{w_j}$ in $D^{(1)}_{\beta(j)j}-D^{(2)}_{\beta(j)j}$ is $\pm1$;
\item every other term $n_w\gamma_w$ of $D^{(1)}_{ij}$ or $D^{(2)}_{ij}$ satisfies $a(w)-a(w_j)\ge u_i-u_{\beta(j)}+1$;
\item $a(w_j)+u_i-u_{\beta(j)}\le41$ for all $i,j$, and $a^P_k\ge16$ for all $k\ne8$.
\end{enumerate}
\end{proposition}

\begin{corollary}\label{prop:cap}
The cap $C$ is acyclic.
\end{corollary}

\begin{proof}
Since $41<3\cdot16=48$, this follows from Propositions~\ref{lem:dominance} and~\ref{prop:capcert}.
\end{proof}

\begin{corollary}\label{prop:m0}
There is an integer $b_0$ such that $M_0=(P^7,\emptyset)$ is acyclic for every $a^P_8\ge b_0$.
\end{corollary}

\begin{proof}
The colors of $I\setminus\{8\}$ span $V$, so Corollary~\ref{prop:cap} and Lemma~\ref{lem:restore} apply.
\end{proof}

Fix $a^P_8\ge b_0$ henceforth.

For the core $X=(Q,\Ip)$ with state $\state_Q$ and fixed-facet weights of Table~\ref{tab:Qweights}, tile labels lie in $U$ and $n=214$ (Proposition~\ref{prop:links}(2)). Table~\ref{tab:basis} specifies $\mathbf e'_j-\mathbf e_j$ for $29$ indices $j$ (with $\mathbf e'_j=\mathbf e_j$ for the other $185$), whose dependency relation $\prec$ is acyclic. Table~\ref{tab:corepivots} lists a bijection $\beta$ of $\{1,\dots,214\}$, paths $w_j$ of length at most $3$ from $x_j$ to $y_{\beta(j)}$, and integers $u_1,\dots,u_{214}$. Exact evaluation (Appendix~\ref{app:computations}, where $1166$, $10947$, and $13074$ nonzero terms in $D^{(1)}$, $D^{(2)}$, and $D'^{\,\le2}$) verifies:

\begin{proposition}\label{prop:corecert}
For the core $X$ (with $(\mathbf e'_j)$ from Table~\ref{tab:basis} and the retractions of Appendix~\ref{app:computations}):
\begin{enumerate}
\item for every $j$, the coefficient of $\gamma_{w_j}$ in $D'^{\,\le2}_{\beta(j)j}$ is $\pm1$;
\item every other term $n\,\gamma_q$ of $D'^{\,\le2}_{ij}$ satisfies $a(q)-a(w_j)\ge u_i-u_{\beta(j)}+1$;
\item for all $i,j$ and every path $p$ occurring in $\mathbf e'_j-\mathbf e_j$ (or the empty path from $x_j$), $a(w_j)-a(p)+u_i-u_{\beta(j)}\le50$, and $a_F\ge17$ for every fixed facet $F$.
\end{enumerate}
\end{proposition}

\begin{theorem}\label{thm:core}
The core $X$ is acyclic.
\end{theorem}

\begin{proof}
Since $50<3\cdot17=51$, this follows from Propositions~\ref{lem:dominance} and~\ref{prop:corecert}.
\end{proof}

\subsection{\texorpdfstring{Duality and acyclicity of the bridge $B$ and end block $E$}{Duality and acyclicity of the bridge B and end block E}}
Let $J=\{i\in\Ip:(\state_Q)_i=1\}$ and $J'=\Ip\setminus J$ (each has $28$ elements). Write $\partial_J X$ and $\partial_{J'}X$ for the unions of the images of $wE_i$ in $\partial X$ for $i\in J$ and $i\in J'$, respectively. Since tiles of $X$ have labels in $U$, the incoming boundary is $\partial_J X$ for $\state_Q$ and $\partial_{J'}X$ for $\state_Q^\sharp$, while $f|_{W^{\Ip}}$ and $\varphi_X$ are the same for both states.

\begin{proposition}\label{prop:dual}
The block $(Q,\Ip)$ with state $\state_Q^\sharp$ and weights of Table~\ref{tab:Qweights} is acyclic.
\end{proposition}

\begin{proof}
Let $R^\pm=\Nov{\Gamma_X}{\pm\varphi_X}$, and write $H_*(X,Y;R^\pm)$ for the homology of $R^\pm\otimes_{\Z\Gamma_X}C_*(\widetilde X,\widetilde Y)$. By Theorem~\ref{thm:core} and Lemma~\ref{lem:model}, we know that $H_*(X,\partial_J X;R^+)=0$, and we must prove $H_*(X,\partial_{J'}X;R^+)=0$.

Because the fixed-facet colors span $U$ and $\one+\eps\in U$, there is $\tau\in W^{\Ip}$ with $\rho(\tau)=\one+\eps$. Since $c_i\cdot(\one+\eps)=1$ for all $i\in\If$, Lemma~\ref{lem:cocycle} gives $\varphi_X(\tau\gamma\tau^{-1})=-\varphi_X(\gamma)$. Left multiplication by $\tau$ on $\widetilde X$ preserves $\widetilde{\partial_J X}$ and intertwines $\gamma$ with $\tau\gamma\tau^{-1}$. Conjugation by $\tau$ gives a ring isomorphism $R^-\to R^+$, so $x\otimes c\mapsto(\tau x\tau^{-1})\otimes\tau c$ is a chain isomorphism from $R^-\otimes C_*(\widetilde X,\widetilde{\partial_J X})$ to $R^+\otimes C_*(\widetilde X,\widetilde{\partial_J X})$. It follows that $H_*(X,\partial_J X;R^-)=0$.

Now remove disjoint open horoball cusp neighbourhoods from $X$ to obtain a compact core $\overline X\subset X$, whose boundary decomposes as $\partial_J\overline X\cup\partial_{J'}\overline X\cup\mathcal E$ with $\mathcal E$ lying on the horospheres. At each ideal vertex of $Q$, cutting along the paired facets turns the Euclidean $6$-cube link into a product of an interval with a $5$-cube (Lemma~\ref{lem:Q}(4)). Each component of $\mathcal E$ is therefore $[0,1]\times T^5$ with ends in $\partial_J\overline X\cup\partial_{J'}\overline X$, and its fundamental group $\Lambda<\Gamma_X$ contains $(r_\alpha r_\beta)^2$ for the five tangent fixed pairs at that vertex. By Lemma~\ref{lem:peripheral}(1) and \eqref{eq:period}, there exists $t\in\Lambda$ with $\varphi_X(t)\ne0$. Since $1-t$ is invertible in $R^-$ and kills $\Z$, we have $\operatorname{Tor}_*^{\Z\Lambda}(R^-,\Z)=0$, and hence $H_*(\mathcal E,\mathcal E\cap\partial_J\overline X;R^-)=0$. Excision and the exact sequence of the triple $(\overline X,\partial_J\overline X\cup\mathcal E,\partial_J\overline X)$ now give
\[
H_*(\overline X,\partial_J\overline X\cup\mathcal E;R^-)\cong H_*(\overline X,\partial_J\overline X;R^-)\cong H_*(X,\partial_J X;R^-)=0.
\]

The compact orientable $7$-manifold $\overline X$ has boundary split along a $5$-manifold into $\partial_J\overline X\cup\mathcal E$ and $\partial_{J'}\overline X$. By Poincar\'e--Lefschetz duality over $\Z\Gamma_X$ \cite[Chapter~2]{Wall99}, $C_*(\widetilde{\overline X},\widetilde{\partial_{J'}\overline X})$ is chain homotopy equivalent to the $\Z\Gamma_X$-dual of $C_*(\widetilde{\overline X},\widetilde{\partial_J\overline X\cup\mathcal E})$, with degrees $k$ and $7-k$ exchanged via the involution $g\mapsto g^{-1}$. This involution extends to an anti-isomorphism $R^-\to R^+$. Since $R^-\otimes C_*(\widetilde{\overline X},\widetilde{\partial_J\overline X\cup\mathcal E})$ is contractible, so is its dual, and therefore $H_*(X,\partial_{J'}X;R^+)\cong H_*(\overline X,\partial_{J'}\overline X;R^+)=0$.
\end{proof}

\begin{proposition}\label{prop:bridge}
There is an integer $L_0$ such that the bridge $B=(Q,\{0,181\})$ is acyclic for every $L\ge L_0$.
\end{proposition}

\begin{proof}
Fix $p\in\Ip\setminus\{0,181\}$. Choose $w_v\in W^{\Ip}$ for $v\in U$, and set $w_v=r_p w_{v+\eps}$ for $v\notin U$. Decompose the differential of $C=C(Q,\{0,181\})$ as $d_0+d_1^{\mathrm f}+\delta$, where $d_1^{\mathrm f}$ consists of the lower faces in fixed directions and $\delta$ of those in the $54$ non-port paired directions. Crossing a fixed facet preserves $v\cdot\eps$, so $d_0+d_1^{\mathrm f}$ is a differential preserving the splitting $C=C_0\oplus C_1$ into labels $v\in U$ and $v\notin U$.

Then $(C_0,d_0+d_1^{\mathrm f})=\Z\Gamma_B\otimes_{\Z\Gamma_X}C(Q,\Ip)$. When $v\notin U$ we have $D(v)=D^\sharp(v+\eps)$, and sending $e_{v,S}$ to $e_{v+\eps,S}$ identifies $(C_1,d_0+d_1^{\mathrm f})$ with $\Z\Gamma_B\otimes_{\Z\Gamma_X}C^\sharp(Q,\Ip)$ via $\gamma\mapsto r_p\gamma r_p$. Because fixed-facet increments are the same at $v$ and $v+\eps$, we have $\varphi_B(\gamma)=\varphi_B(r_p\gamma r_p)=\varphi_X(\gamma)$ for every $\gamma\in\Gamma_X$. Both homomorphisms thus extend to valuation-preserving maps from $\Nov{\Gamma_X}{\varphi_X}$ into $R=\Nov{\Gamma_B}{\varphi_B}$. Contractions $\Theta$ and $\Theta^\sharp$ of the two core complexes (Theorem~\ref{thm:core} and Proposition~\ref{prop:dual}) then induce a contraction $\Theta_0$ of $(R\otimes C,d_0+d_1^{\mathrm f})$. On $C_0$ the basis heights $\eta_v$ match those of $X$, whereas on $C_1$ they differ by the constant $f(r_p)=\pm L$. Consequently $\|\Theta_0\|=\max\{\|\Theta\|,\|\Theta^\sharp\|\}$ does not depend on $L$. Since $\|\delta\|=-L$, Lemma~\ref{lem:perturb} shows that $B$ is acyclic whenever $L>\|\Theta_0\|$.
\end{proof}

We now fix $L\ge L_0$.

\begin{corollary}\label{prop:end}
There is an integer $a_0$ such that the end block $E=(Q,\{0\})$ is acyclic for every $a_*\ge a_0$.
\end{corollary}

\begin{proof}
The colors of $I\setminus\{0,181\}$ already span $V$, so Proposition~\ref{prop:bridge} and Lemma~\ref{lem:restore} apply.
\end{proof}

We also fix $a_*\ge a_0$.

\begin{remark}\label{rem:even}
The same argument shows that $(Q,\{181\})$ is acyclic when the weight of $E_0$ is large enough. This block appears at the end of $M_m$ for even $m\ge2$, so Section~\ref{sec:gluing} applies to even $m$ as well.
\end{remark}

\section{\texorpdfstring{Gluing the blocks and proof of Theorem~\ref{thm:main}}{Gluing the blocks and proof of Theorem A}}\label{sec:gluing}

Let $m=0$ or let $m\ge1$ be odd, and let $\alpha\colon\Gamma_m\to\Z$ be the character of Section~\ref{sec:data} with the parameters $a^P_8,L,a_*$ fixed in Section~\ref{sec:novikov}. Write $p\colon\Hy^7\to M_m$ for the covering projection.

\subsection{The block decomposition and the graph character}
Write $Q_0=P^7$, $N_k=p(W_m\cdot Q_k)$ for $0\le k\le m$, and $Z_k=p(W_m\cdot(Q_{k-1}\cap Q_k))$ for $1\le k\le m$. Then $M_m$ is the union of $N_0,\dots,N_m$, with $N_{k-1}\cap N_k=Z_k$ and $N_k\cap N_l=\emptyset$ whenever $|k-l|\ge2$. Any facet meeting an inlet has color in $U$, so $v\cdot\eps$ is constant on each component of $Z_k$. Denote by $Z_k^\epsilon$ the union of inlets of tiles with $v\cdot\eps=\epsilon$, and regard $G_k=\pi_1(N_k)$ and $G_k^\epsilon=\pi_1(Z_k^\epsilon)$ as subgroups of $\Gamma_m$ defined up to conjugacy.

\begin{lemma}\label{lem:blocks}
\begin{enumerate}
\item Each $N_k$ is a connected submanifold with totally geodesic boundary $Z_k\sqcup Z_{k+1}$ (taking $Z_0=Z_{m+1}=\emptyset$), and each $Z_k^\epsilon$ is connected.
\item With $\partial^-N_k=Z_k^0\sqcup Z_{k+1}^1$, we have $H_*(N_k,\partial^-N_k;\Nov{G_k}{\alpha})=0$.
\item The restrictions $\alpha|_{G_k}$ and $\alpha|_{G_k^\epsilon}$ are nonzero with finitely presented kernels, and $\alpha(G_0)=\Z$.
\end{enumerate}
\end{lemma}

\begin{proof}
When $m=0$ we have $N_0=M_0$, empty boundary, and $G_0=\Gamma_0$, so everything follows from Corollary~\ref{prop:m0}, Lemma~\ref{lem:model}, Corollary~\ref{cor:fp}, and Lemma~\ref{lem:peripheral}(3).

Suppose $m\ge1$. Set $\Pi_0=\{8\}$ with $h_0=1$, $\Pi_k=\{0,181\}$ for $1\le k<m$, and $\Pi_m=\{0\}$ (since $m$ is odd, the inlet of $Q_m$ is $h_m(E_0)$). Let $W^{(k)}<W_m$ be the subgroup generated by reflections in facets of $S_m$ containing non-port facets of $Q_k$. By Lemma~\ref{lem:Sm}, we have $W^{(k)}=h_k W^{\Pi_k}h_k^{-1}$ with matching colors, and hence $\Gamma_m\cap W^{(k)}=h_k\Gamma^{\Pi_k}h_k^{-1}$. Since the non-port colors span $V$, the map $x\mapsto p(h_k x)$ induces an isometry from $N(Q_k,\Pi_k)$ (that is, $C$, $B$, or $E$) onto $N_k$. This identifies $\alpha|_{G_k}$ with $\varphi_C$ ($k=0$), $\varphi_B$ or $\varphi_E$ ($k$ odd), or $\varphi_B^\sharp$ ($k$ even).

By Lemma~\ref{lem:compatible}, the inlet of $Q_k$ ($h_k(E_0)$ for $k$ odd, $h_k(E_{181})$ for $k$ even) has state $1$ and color $\eps$, so it is incoming if and only if $v\cdot\eps=0$. The outlet of $Q_k$ ($k<m$) and $F_8$ of $P^7$ have state $0$ and color $\eps$, so they are incoming if and only if $v\cdot\eps=1$. Hence the incoming boundary of $N_k$ is $\partial^-N_k=Z_k^0\sqcup Z_{k+1}^1$.

(1) Each block is connected. Its boundary components over a port $i$ correspond to $V/\rho(W_{(i)})$, which has $2$ elements (distinguished by $v\cdot\eps$) because the colors of facets meeting any port span $U$ (Appendix~\ref{app:computations}).

(2) This follows from Lemma~\ref{lem:model}, Corollaries~\ref{prop:cap} and~\ref{prop:end}, Proposition~\ref{prop:bridge}, and Lemma~\ref{lem:sharp} (for even $k$).

(3) The kernels of $\alpha|_{G_k}$ and $\alpha|_{G_k^\epsilon}$ are finitely presented by Corollary~\ref{cor:fp} and Lemma~\ref{lem:sharp}, since facets meeting an inlet are fixed and have the same states in $\state_Q$ and $\state_Q^\sharp$. Every inlet $E_i$ contains an ideal vertex whose five tangent fixed pairs all meet $E_i$, so $\alpha|_{G_k^\epsilon}\ne0$ by Lemma~\ref{lem:peripheral}(1) and \eqref{eq:period}. It follows that $\alpha|_{G_k}\ne0$ as well. Finally, $\alpha(G_0)=\varphi_C(\Gamma_C)=\Z$ by Lemma~\ref{lem:peripheral}(3).
\end{proof}

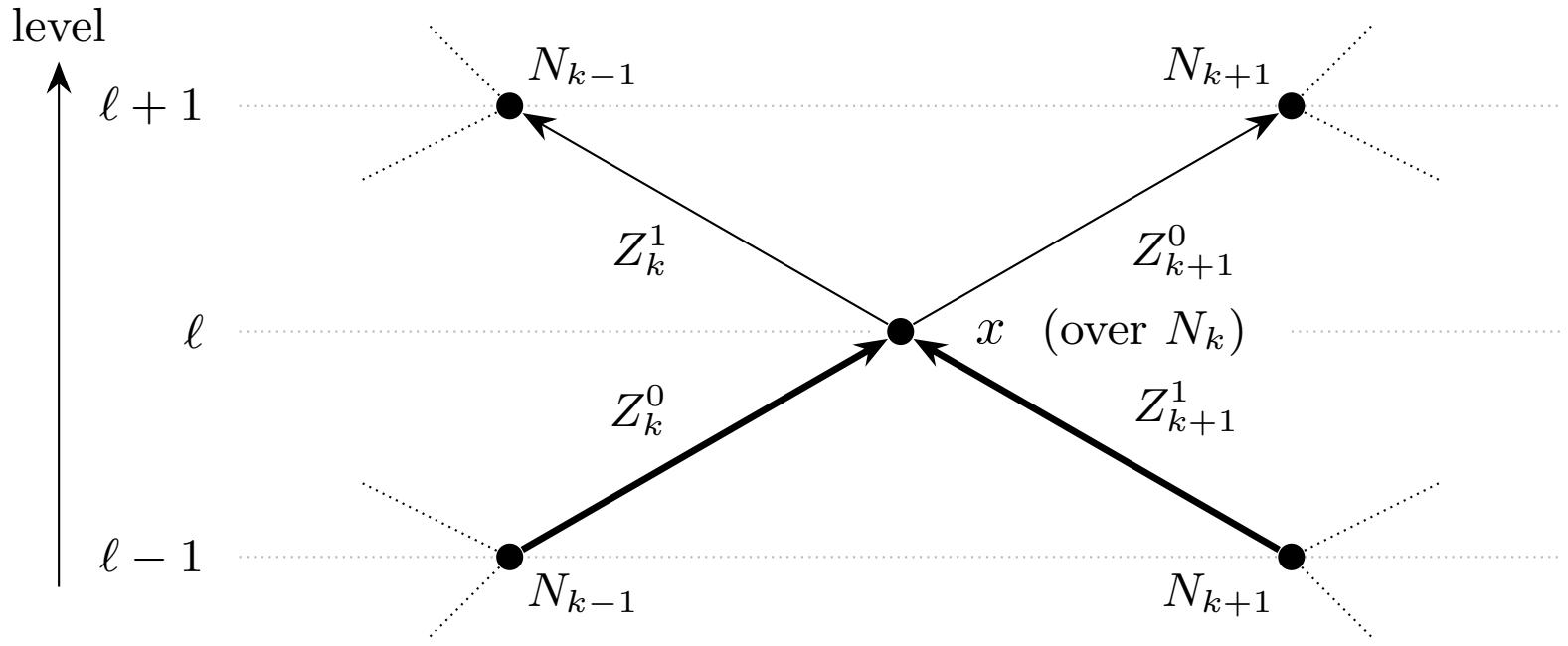
\begin{figure}[ht]
\centering
\begin{tikzpicture}[font=\small,x=1cm,y=1cm,>={Stealth[length=2.2mm]},vx/.style={circle,fill=black,inner sep=1.8pt}]
\foreach \y/\t in {0/{\ell-1},3/{\ell+1}} {\draw[gray!55,densely dotted] (-4.4,\y)--(4.4,\y); \node[font=\footnotesize,anchor=east] at (-4.5,\y) {$\t$};}
\draw[gray!55,densely dotted] (-4.4,1.5)--(-0.2,1.5); \draw[gray!55,densely dotted] (2.6,1.5)--(4.4,1.5); \node[font=\footnotesize,anchor=east] at (-4.5,1.5) {$\ell$};
\draw[->] (-5.6,-0.2)--(-5.6,3.3) node[above,font=\footnotesize] {level};
\node[vx] (x) at (0,1.5) {};
\node[vx] (a) at (-2.6,0) {}; \node[vx] (b) at (2.6,0) {};
\node[vx] (c) at (-2.6,3) {}; \node[vx] (d) at (2.6,3) {};
\draw[line width=1.3pt,->] (a)--(x) node[pos=0.45,above left=-1pt,font=\footnotesize] {$Z_k^0$};
\draw[line width=1.3pt,->] (b)--(x) node[pos=0.45,above right=-1pt,font=\footnotesize] {$Z_{k+1}^1$};
\draw[->] (x)--(c) node[pos=0.55,below left=-1pt,font=\footnotesize] {$Z_k^1$};
\draw[->] (x)--(d) node[pos=0.55,below right=-1pt,font=\footnotesize] {$Z_{k+1}^0$};
\node[right=13pt,inner sep=1pt] at (x) {$x$ \ {\footnotesize(over $N_k$)}};
\node[anchor=north west,font=\footnotesize,inner sep=3pt] at (a) {$N_{k-1}$}; \node[anchor=north east,font=\footnotesize,inner sep=3pt] at (b) {$N_{k+1}$};
\node[anchor=south west,font=\footnotesize,inner sep=3pt] at (c) {$N_{k-1}$}; \node[anchor=south east,font=\footnotesize,inner sep=3pt] at (d) {$N_{k+1}$};
\foreach \p/\s in {a/-1,b/1} {\draw[densely dotted] (\p)--++(0.55*\s,-0.55); \draw[densely dotted] (\p)--++(1.0*\s,0.5);}
\foreach \p/\s in {c/-1,d/1} {\draw[densely dotted] (\p)--++(0.55*\s,0.55); \draw[densely dotted] (\p)--++(1.0*\s,-0.5);}
\end{tikzpicture}
\caption{A part of the Bass--Serre tree $\cT$ around a vertex $x$ over $N_k$ \mbox{($0<k<m$)}, with one edge of each kind drawn and the level increasing upwards. The edges that end at $x$ (thick) are the components of the incoming boundary of the corresponding copy of $\widetilde N_k$.}\label{fig:tree}
\end{figure}

Let $\cG$ be the oriented graph with vertices $0,\dots,m$ and edges $e_k^0\colon(k-1)\to k$ and $e_k^1\colon k\to(k-1)$ ($1\le k\le m$). The vertex $k$ stands for $N_k$ and $e_k^\epsilon$ for $Z_k^\epsilon$, so every edge points into the block for which $Z_k^\epsilon$ is incoming (Figure~\ref{fig:chain}). Let $\cT$ be the Bass--Serre tree dual to the decomposition of $\Hy^7$ along the preimages of the $Z_k$, oriented compatibly with $\cT/\Gamma_m=\cG$. Choose a \emph{level} function $\ell\colon V(\cT)\to\Z$ that increases by $1$ along every oriented edge (Figure~\ref{fig:tree}). This defines a homomorphism
\[
\psi\colon\Gamma_m\to\Z,\qquad\psi(g)=\ell(gx)-\ell(x),
\]
which vanishes on all vertex and edge stabilizers (conjugates of $G_k$ and $G_k^\epsilon$, with $\psi=0$ when $m=0$). For $N\in\Z$, put $\alpha_N=\alpha+N\psi\colon\Gamma_m\to\Z$.

\subsection{\texorpdfstring{Novikov acyclicity of $M_m$}{Novikov acyclicity of Mm}}

\begin{theorem}\label{thm:glue}
\begin{enumerate}
\item For every $N\in\Z$, $\alpha_N$ is surjective and its restriction to each cusp group of $M_m$ is nonzero.
\item For every $N\in\Z$, $\ker\alpha_N$ is finitely presented.
\item There exists $N_0$ such that for all $N\ge N_0$ and every finite CW complex $Y\simeq M_m$,
\[
H_*\bigl(\Nov{\Gamma_m}{\alpha_N}\otimes_{\Z\Gamma_m}C_*(\widetilde Y)\bigr)=0=H_*\bigl(\Nov{\Gamma_m}{-\alpha_N}\otimes_{\Z\Gamma_m}C_*(\widetilde Y)\bigr).
\]
\end{enumerate}
\end{theorem}

\begin{proof}
(1) Since $\psi$ vanishes on $G_0$, we have $\alpha_N(G_0)=\alpha(G_0)=\Z$ by Lemma~\ref{lem:blocks}(3). A cusp group $\Lambda$ is the stabilizer in $\Gamma_m$ of $g\xi$ for some $g\in W_m$ and ideal vertex $\xi$ of a piece $Q_k$. At least five tangent pairs $(F_\alpha,F_\beta)$ at $\xi$ come from non-port facets of $Q_k$. For any such pair, the element $t=g(r_\alpha r_\beta)^2g^{-1}\in\Lambda$ stabilizes the vertex $g\cdot\widetilde N_k$ of $\cT$. Hence $\psi(t)=0$, and \eqref{eq:period} gives $\alpha_N(t)=\alpha(t)=2(\alpha(\rho_m(g),\alpha)-\alpha(\rho_m(g),\beta))$. Because Lemma~\ref{lem:peripheral}(1)--(2) holds for every $v\in V$, it applies to both $\state_Q$ and $\state_Q^\sharp$ and guarantees that this period is nonzero for at least one pair.

(2) Let $K=\ker\alpha_N$ act on the tree $\cT$. Vertex stabilizers are isomorphic to $\ker(\alpha|_{G_k})$ and edge stabilizers to $\ker(\alpha|_{G_k^\epsilon})$, all finitely presented by Lemma~\ref{lem:blocks}(3). Because $\alpha(G_k)$ and $\alpha(G_k^\epsilon)$ are nonzero subgroups of $\alpha_N(\Gamma_m)=\Z$, they have finite index in $\Z$, so there are finitely many $K$-orbits of vertices and edges. Therefore $K$ is the fundamental group of a finite graph of finitely presented groups, and is itself finitely presented.

(3) Choose disjoint $W_m$-invariant open horoball cusp neighbourhoods, and let $\Mbar\subset M_m$ be their complement, with a finite triangulation in which all $N_k\cap\Mbar$ and $Z_k^\epsilon\cap\Mbar$ are subcomplexes. By homotopy invariance of Novikov homology, it suffices to take $Y=\Mbar$. For each vertex $x$ and edge $e$ of $\cT$, let $C(x)$ and $C(e)$ be the cellular chain complexes of the corresponding subspaces of $\widetilde Y$. The exact sequence
\[
0\longrightarrow\bigoplus_e C(e)\xrightarrow{\ \Phi^{\mathrm t}-\Phi^{\mathrm s}\ }\bigoplus_x C(x)\longrightarrow C_*(\widetilde Y)\longrightarrow0
\]
(where $\Phi^{\mathrm t}$ and $\Phi^{\mathrm s}$ are the inclusions of $C(e)$ into its target and source) shows that $C_*(\widetilde Y)$ is chain homotopy equivalent to $\operatorname{Cone}(\Phi^{\mathrm t}-\Phi^{\mathrm s})$. Write its differential as $d^{\mathrm t}+\delta$, where $d^{\mathrm t}$ is the differential of $\operatorname{Cone}(\Phi^{\mathrm t})$ and $\delta$ comes from $-\Phi^{\mathrm s}$.

Put $R=\Nov{\Gamma_m}{\alpha_N}$. The complex $\operatorname{Cone}(\Phi^{\mathrm t})$ splits over $x\in V(\cT)$ as the direct sum of the mapping cones of $\bigoplus_{\mathrm{t}(e)=x}C(e)\hookrightarrow C(x)$, each chain homotopy equivalent to the relative complex of $x$ modulo its incoming boundary. Choose one vertex $x_k$ over each $N_k$ with stabilizer $G_k$. By Lemma~\ref{lem:blocks}(2), each summand at $x_k$ has a contraction $\Theta_k$ over $\Nov{G_k}{\alpha}\subseteq R$, and these induce a contraction $\Theta$ of $(R\otimes\operatorname{Cone}(\Phi^{\mathrm t}),d^{\mathrm t})$.

Choose a $\Z\Gamma_m$-basis of $\operatorname{Cone}(\Phi^{\mathrm t}-\Phi^{\mathrm s})$ consisting of cells in the vertices $x_k$ and in edges with target $x_k$ or $x_{k-1}$. Define the level $\ell(b)$ of a basis cell $b$ as the level of its vertex (or target vertex, for an edge cell), and give $b$ the height $N\ell(b)$. Using $\ell(gb')=\ell(b')+\psi(g)$, a term $g\,b'$ in the image of $b$ shifts the height by
\[
\alpha_N(g)+N\ell(b')-N\ell(b)=\alpha(g)+N\bigl(\ell(gb')-\ell(b)\bigr).
\]
Because $d^{\mathrm t}$ and $\Theta$ keep cells at the same level, $\|\Theta\|$ is independent of $N$. By contrast, $\delta$ maps each edge cell from its target to its source, dropping the level by $1$, so $\|\delta\|=c-N$ for a constant $c$ not depending on $N$. Lemma~\ref{lem:perturb} then implies that $\operatorname{Cone}(\Phi^{\mathrm t}-\Phi^{\mathrm s})$ is $R$-acyclic whenever $N>c+\|\Theta\|$.

For $-\alpha_N$, pick $\tau\in W_m$ with $\rho_m(\tau)=\one$. Lemma~\ref{lem:cocycle} gives $\alpha(\tau\gamma\tau^{-1})=-\alpha(\gamma)$. Left translation by $\tau$ preserves the horoball truncation and each piece type $Q_k$, while exchanging $Z_k^0$ and $Z_k^1$ because $\one\cdot\eps=1$. It therefore reverses the orientation of every edge of $\cT$, whence $\psi(\tau\gamma\tau^{-1})=-\psi(\gamma)$ and $\alpha_N(\tau\gamma\tau^{-1})=-\alpha_N(\gamma)$. Conjugating coefficients by $\tau$ and acting by $\tau$ on $\widetilde Y$ yields an isomorphism between $H_*(\Nov{\Gamma_m}{-\alpha_N}\otimes C_*(\widetilde Y))$ and $H_*(\Nov{\Gamma_m}{\alpha_N}\otimes C_*(\widetilde Y))=0$.
\end{proof}

\subsection{\texorpdfstring{Proof of Theorem~\ref{thm:main}}{Proof of Theorem A}}\label{sec:proof}
Fix $N\ge N_0$ as in Theorem~\ref{thm:glue} and set $\varphi=\alpha_N$. Let $\Mbar\subset M_m$ be the compact core, a compact smooth $7$-manifold with boundary $\partial\Mbar=\bigsqcup_c T_c$ ($T_c\cong T^6$, $252+1458m$ components) and interior diffeomorphic to $M_m$, and let $\Mbar_\varphi$ be its infinite cyclic cover corresponding to $\ker\varphi$.

Recall that a space $Y$ is \emph{finitely dominated} if there exist a finite CW complex $Z$ and maps $i\colon Y\to Z$ and $r\colon Z\to Y$ such that $r\circ i$ is homotopic to the identity. We use Farrell's fibering theorem \cite{Farrell71} (see also \cite{Farrell72}), which does not require the boundary to be connected.

\begin{theorem}\label{thm:farrell}
Let $M$ be a compact connected smooth manifold of dimension at least $6$, and let $g\colon M\to S^1$ be a smooth map that induces a surjection on fundamental groups and restricts to a smooth fiber bundle on $\partial M$. Suppose that the infinite cyclic cover of $M$ determined by $\ker g_*$ is finitely dominated and that $\Wh(\pi_1 M)=0$. Then $g$ is homotopic relative to $\partial M$ to a smooth fiber bundle.
\end{theorem}

Finite domination of an infinite cyclic cover is characterized in terms of Novikov homology by Ranicki \cite{Ranicki95} and Sch\"utz \cite[Theorem~4.9]{Schuetz06}.

\begin{theorem}\label{thm:schuetz}
Let $Y$ be a finite connected CW complex with fundamental group $G$ and universal cover $\widetilde Y$. Let $\chi\colon G\to\Z$ be an epimorphism, and let $Y_\chi$ be the infinite cyclic cover of $Y$ with fundamental group $\ker\chi$. Then $Y_\chi$ is finitely dominated if and only if $\ker\chi$ is finitely presented and
\[
H_*\bigl(\Nov{G}{\chi}\otimes_{\Z G}C_*(\widetilde Y)\bigr)=0=H_*\bigl(\Nov{G}{-\chi}\otimes_{\Z G}C_*(\widetilde Y)\bigr).
\]
\end{theorem}

\begin{lemma}\label{lem:boundary}
There is a smooth map $f_1\colon\Mbar\to S^1$ inducing $\varphi$ on $\pi_1$ whose restriction to each boundary torus $T_c$ is a smooth fiber bundle with fiber consisting of $d_c$ copies of $T^5$, where $d_c=[\Z:\varphi(\pi_1T_c)]\ge1$.
\end{lemma}

\begin{proof}
Write $T_c=\R^6/\Lambda_c$ ($\Lambda_c\cong\Z^6$). By Theorem~\ref{thm:glue}(1), $\varphi(\Lambda_c)=d_c\Z$ with $d_c\ge1$. The linear extension $\lambda_c\colon\R^6\to\R$ of $\varphi|_{\Lambda_c}$ descends to a smooth fiber bundle $\bar\lambda_c\colon T_c\to\R/\Z=S^1$ with fiber $\lambda_c^{-1}(\Z)/\Lambda_c\cong\bigsqcup_{d_c}T^5$. Starting from any smooth map $f_0\colon\Mbar\to S^1$ inducing $\varphi$, we have $f_0|_{T_c}\simeq\bar\lambda_c$ for each $c$, and a collar homotopy on $\partial\Mbar$ deforms $f_0$ to a smooth map $f_1$ with $f_1|_{T_c}=\bar\lambda_c$.
\end{proof}

\begin{proof}[Proof of Theorem~\ref{thm:main}]
Parts~(1) and~(4) are Proposition~\ref{prop:cover} (with Lemma~\ref{lem:Sm}) and Theorem~\ref{prop:commens}.

For~(2), we apply Theorem~\ref{thm:farrell} to $\Mbar$ ($\dim\Mbar=7\ge6$) and the map $f_1\colon\Mbar\to S^1$ of Lemma~\ref{lem:boundary}. By Theorem~\ref{thm:glue}, $\varphi$ is surjective with finitely presented kernel and vanishing Novikov homology for both $\pm\varphi$, so $\Mbar_\varphi$ is finitely dominated by Theorem~\ref{thm:schuetz}. Moreover, $\Gamma_m$ is torsion-free (Proposition~\ref{prop:cover}(1)) and acts freely, cocompactly and isometrically on the $7$-dimensional $\mathrm{CAT}(0)$ Davis complex of $W_m$, so $\Wh(\Gamma_m)=0$ by \cite[Theorem~B]{BL12}. Theorem~\ref{thm:farrell} therefore homotopes $f_1$ relative to $\partial\Mbar$ to a smooth fiber bundle $\Mbar\to S^1$ with compact $6$-manifold fiber $F$. Restricting to $\operatorname{int}\Mbar\cong M_m$ yields a smooth fiber bundle $M_m\to S^1$ with fiber $\operatorname{int}F$.

For~(3), the homotopy exact sequence shows that $F$ is connected ($\varphi$ is surjective) and aspherical ($\Mbar$ is aspherical). It is orientable since $\Mbar$ and $S^1$ are, and $\partial F=F\cap\partial\Mbar$ consists of $\sum_c d_c$ copies of $T^5$ by Lemma~\ref{lem:boundary}.
\end{proof}

\begin{remark}
The fibration represents $\alpha+N\psi$, where $a^P_8,L,a_*,N$ are bounded in terms of the norms of the contractions involved. We have not computed explicit values, so $\sum_c d_c$ is not determined here.
\end{remark}

\section{Further remarks and open questions}\label{sec:questions}

The hybrid construction of Section~\ref{sec:manifold} has a $5$-dimensional analogue. Let $P^5\subset\Hy^5$ and $P^6\subset\Hy^6$ be the right-angled polytopes of \cite{PV05,ERT12} ($16$ and $27$ facets), with normals obtained by permuting the first $n$ coordinates of $(-1,0,\dots,0;0)$, $(1,1,0,\dots,0;1)$, and $(1,1,1,1,1,0,\dots,0;2)$. The section $Q^5=P^6\cap\{x_1=x_2\}\subset\Hy^5$ has $21$ facets. Its facet $E$ with normal $(-1,0,0,0,0;0)$ is isometric to $P^4$ and matches a facet of $P^5$ (after reflecting $P^5$ in $\{x_1=0\}$), while the facet $H$ with normal $(\sqrt2,0,0,0,0;1)$ is ultraparallel to $E$ (inner product $-\sqrt2$). Doubling along $H$ gives a $28$-facet right-angled polytope $B^5=Q^5\cup r_H(Q^5)$ with ultraparallel ports $E,r_H(E)$, and gluing $P^5$ to a chain of $m$ copies of $B^5$ yields a right-angled polytope $L_m\subset\Hy^5$.

\begin{proposition}\label{prop:dim5}
For every $m\ge1$, $L_m$ is a finite-volume right-angled polytope in $\Hy^5$ with $16+16m$ facets. For $m\ne m'$, the reflection groups of $L_m$ and $L_{m'}$ are not commensurable in $\Isom(\Hy^5)$ (even up to conjugacy). In particular, reflection groups of finite-volume right-angled polytopes in $\Hy^5$ fall into infinitely many commensurability classes.
\end{proposition}

\begin{proof}
As in Lemmas~\ref{lem:Q} and~\ref{lem:Sm}, $L_m$ is right-angled of finite volume. Each port of $B^5$ meets $10$ facets, so each copy of $B^5$ contributes $28-2-10=16$ facets. By the argument of Lemma~\ref{lem:arith}(1), the reflection groups of $P^5$ and $Q^5$ are commensurable with $O(q,\Z)$ for $x_1^2+\dots+x_5^2-t^2$ and $2x_1^2+x_2^2+\dots+x_5^2-t^2$. These $6$-ary forms have discriminants $-1$ and $-2$ in $\Q^\times/(\Q^\times)^2$, so they are not similar over $\Q$. Suppose a torsion-free lattice $\Lambda$ has finite indices $d,d'$ in conjugates of the reflection groups of $L_m$ and $L_{m'}$. In $\Hy^5/\Lambda$, the $P^5$- and $B^5$-pieces double along their boundaries to finite covers of the reflection orbifolds of $P^5$ and $Q^5$. Applying \cite[Lemma~3.5]{GL14} as in Theorem~\ref{prop:commens} then yields $d=d'$ and $m=m'$.
\end{proof}

Our arguments do not answer whether $L_m$ has a fibered finite cover. Farrell's fibering theorem applies only in dimension at least $6$, while the direct construction of \cite{IMM23} needs every ascending and descending link to be collapsible and cannot use the perturbation methods of Sections~\ref{sec:novikov} and~\ref{sec:gluing}. In dimensions $n\ge9$, finite-volume right-angled hyperbolic polytopes do not exist \cite{PV05,Dufour10,Kolpakov26}, so one would need to start from non-right-angled Coxeter polytopes or from lattices not commensurable with reflection groups.

\begin{question}\label{q:dim5}
Are there finite-volume hyperbolic $5$-manifolds that fiber over the circle in infinitely many commensurability classes? For example, is there, for every $m$, a manifold tessellated by finitely many copies of $L_m$ that fibers over the circle?
\end{question}

\begin{question}\label{q:dim}
Does there exist a finite-volume hyperbolic $n$-manifold that fibers over the circle for odd $n\ge9$ (cf.\ \cite[Question~26]{IMM23})?
\end{question}

\section*{Use of computers}
Several finite steps in the proof rely on exact computations over $\Z$, $\Z[\sqrt2]$ and $\F_2$. These cover the polytope incidences, color spans, and compatibility checks (Lemmas~\ref{lem:P7}, \ref{lem:Q}, and~\ref{lem:compatible}), the simple connectivity and homology of the ascending and descending links (Proposition~\ref{prop:links}), the cusp periods and surjectivity of $\varphi_C$ (Lemma~\ref{lem:peripheral}), and the leading-term bounds for the cap $C$ and the core $X$ (Propositions~\ref{prop:capcert} and~\ref{prop:corecert}). All explicit data are listed in Appendix~\ref{app:data}, and the verification procedures are described in Appendix~\ref{app:computations}. The main matrix and leading-term checks were carried out independently by two separate programs included in the ancillary files.

\section*{Use of AI tools}
The proof strategy of this paper is due to the author. ChatGPT (OpenAI) assisted in developing the proofs and writing the code for the computations. In particular, the algebraic-topological computations of Novikov homology in Section~\ref{sec:novikov} were carried out primarily with ChatGPT. Claude (Anthropic) and Gemini (Google) were used to polish the manuscript. The author takes full responsibility for the mathematical content of this paper.

\section*{Acknowledgements}
The author thanks Fangting Zheng for introducing the author to related topics five years ago, and for helpful discussions during the writing of this paper.

\appendix
\section{Data}\label{app:data}

\emph{The polytope $P^7$.} A vector $(x_1,\dots,x_7;t)\in\Z^8$ is written as the word $x_1x_2\cdots x_7\,|\,t$, and a bar denotes a negative entry. For example $\bar1011011\,|\,2$ stands for $(-1,0,1,1,0,1,1;2)$. Table~\ref{tab:normals} lists the normals $e_j$ of the facets of $P^7$. The color $c_j\in V=\F_2^8$ is obtained by reducing the eight entries modulo $2$. An element $v=(v_1,\dots,v_8)$ of $V$ is written as the binary word $v_1v_2\cdots v_8$.

\begin{table}[ht]
\centering\footnotesize\setlength{\tabcolsep}{3.5pt}
\begin{tabular}{c|cccccccc}
\toprule
$n\backslash r$&0&1&2&3&4&5&6&7\\\midrule
0&$0\bar100000\,|\,0$&$00\bar10000\,|\,0$&$\bar1011110\,|\,2$&$\bar1101110\,|\,2$&$0010001\,|\,1$&$0100001\,|\,1$&$\bar1121111\,|\,3$&$\bar1211111\,|\,3$\\
1&$1000000\,|\,0$&$000\bar1000\,|\,0$&$0001001\,|\,1$&$\bar1000001\,|\,1$&$0111110\,|\,2$&$\bar1110110\,|\,2$&$\bar1112111\,|\,3$&$\bar2111111\,|\,3$\\
2&$0010010\,|\,1$&$0010100\,|\,1$&$0100010\,|\,1$&$0100100\,|\,1$&$\bar1011011\,|\,2$&$\bar1011101\,|\,2$&$\bar1101011\,|\,2$&$\bar1101101\,|\,2$\\
3&$0011000\,|\,1$&$0011111\,|\,2$&$\bar1010000\,|\,1$&$\bar1010111\,|\,2$&$0101000\,|\,1$&$0101111\,|\,2$&$\bar1100000\,|\,1$&$\bar1100111\,|\,2$\\
4&$0001010\,|\,1$&$0001100\,|\,1$&$\bar1000010\,|\,1$&$\bar1000100\,|\,1$&$0111011\,|\,2$&$0111101\,|\,2$&$\bar1110011\,|\,2$&$\bar1110101\,|\,2$\\
5&$0000\bar100\,|\,0$&$0000011\,|\,1$&$00000\bar10\,|\,0$&$0000101\,|\,1$&$\bar1111010\,|\,2$&$\bar1111121\,|\,3$&$\bar1111100\,|\,2$&$\bar1111211\,|\,3$\\
6&$000000\bar1\,|\,0$&$0000110\,|\,1$&$\bar1001000\,|\,1$&$\bar1001111\,|\,2$&$0110000\,|\,1$&$0110111\,|\,2$&$\bar1111001\,|\,2$&$\bar1111112\,|\,3$\\\bottomrule
\end{tabular}
\caption{The outward unit normals $e_{8n+r}$ of the facets $F_{8n+r}$ of $P^7$.}\label{tab:normals}
\end{table}

\begin{table}[ht]
\centering\footnotesize\setlength{\tabcolsep}{4pt}
\begin{tabular}{c|cccccccc}
\toprule
$n\backslash r$&0&1&2&3&4&5&6&7\\\midrule
0&27&27&27&27&18&22&27&27\\
1&---&19&27&27&17&27&27&27\\
2&27&27&27&24&17&27&27&27\\
3&18&17&27&27&24&25&27&17\\
4&23&20&27&27&25&23&17&27\\
5&27&19&17&27&27&27&16&27\\
6&25&24&27&17&17&17&27&27\\\bottomrule
\end{tabular}
\caption{The weights $a^P_{8n+r}$ of the facets of $P^7$. The weight of the port $F_8$ is not listed. It is the integer $a^P_8$ of Section~\ref{sec:data}.}\label{tab:capweights}
\end{table}

\emph{The polytope $Q$.} The facets $E_0,\dots,E_{181}$ are numbered as in Section~\ref{sec:manifold}. The state $\state_Q\in\F_2^{182}$ is
\[
\begin{array}{rl}
\state_Q={}&\texttt{11011100 11011110 00101000 01111100 00100000 10000010}\\
&\texttt{11010111 10101011 00101010 00010100 10111110 11111011}\\
&\texttt{11000001 11101011 10000100 11000101 10101001 00011001}\\
&\texttt{01000110 10101001 11010110 01110110 101000},
\end{array}
\]
written in blocks of $8$ bits, the last block having $6$ bits. The bit in position $r$ of block $n$, $0\le r\le7$, $0\le n\le22$, is the coordinate of the facet $E_{8n+r}$. Table~\ref{tab:Qweights} lists the weights of the fixed facets.

\begin{table}[ht]
\centering\footnotesize\setlength{\tabcolsep}{3.5pt}
\begin{tabular}{r|cccccccccccccccccc}
\toprule
1--18&27&27&19&27&17&25&19&27&24&27&23&20&18&27&27&18&17&22\\
19--36&27&24&24&25&17&17&25&23&17&17&23&19&27&20&19&27&21&17\\
37--54&17&20&18&22&19&17&27&17&27&24&17&17&24&17&27&22&22&27\\
55--72&27&27&19&24&23&20&24&19&27&20&19&22&27&26&20&24&27&19\\
73--90&27&25&27&27&22&26&24&19&24&22&17&27&27&25&21&22&20&27\\
91--108&27&17&27&22&17&27&18&22&27&20&27&27&20&24&27&21&27&27\\
109--126&19&17&23&27&27&21&23&21&27&21&27&20&27&27&27&26&27&27\\\bottomrule
\end{tabular}
\caption{The weights $a_i$ of the fixed facets $E_i$ of $Q$, $i=1,\dots,126$. Each row lists eighteen consecutive weights.}\label{tab:Qweights}
\end{table}

\begin{table}[ht]
\centering\small
\begin{tabular}{rcrclr}
\toprule
$j$&$x_j$&$i=\beta(j)$&$y_i$&$w_j$&$u_i$\\\midrule
1&\texttt{10111100}&14&\texttt{00001010}&$r_{20}$&$0$\\
2&\texttt{11011111}&15&\texttt{00010001}&$r_{31}$&$0$\\
3&\texttt{10111111}&13&\texttt{11110110}&$r_{19}$&$0$\\
4&\texttt{10100000}&10&\texttt{00010110}&$r_{20}$&$0$\\
5&\texttt{01010000}&17&\texttt{11010010}&$r_{37}r_{46}$&$0$\\
6&\texttt{00110011}&9&\texttt{11001001}&$r_{46}r_{48}$&$0$\\
7&\texttt{10001001}&22&\texttt{10111000}&$r_{24}$&$0$\\
8&\texttt{11101001}&4&\texttt{01011111}&$r_{20}$&$-3$\\
9&\texttt{00110110}&5&\texttt{11001100}&$r_{46}r_{48}$&$0$\\
10&\texttt{10101010}&20&\texttt{01001011}&$r_{38}r_{41}$&$0$\\
11&\texttt{10001010}&21&\texttt{01000100}&$r_{31}$&$0$\\
12&\texttt{11110101}&3&\texttt{01000011}&$r_{20}$&$0$\\
13&\texttt{00001001}&1&\texttt{01000000}&$r_{19}$&$0$\\
14&\texttt{01101010}&24&\texttt{01011011}&$r_{24}$&$0$\\
15&\texttt{11101110}&2&\texttt{00100000}&$r_{31}$&$0$\\
16&\texttt{10110001}&19&\texttt{00101000}&$r_{41}r_{51}$&$0$\\
17&\texttt{00101101}&6&\texttt{10101111}&$r_{37}r_{46}$&$0$\\
18&\texttt{10010010}&23&\texttt{11011011}&$r_{19}$&$-7$\\
19&\texttt{10110100}&8&\texttt{01010101}&$r_{38}r_{41}$&$0$\\
20&\texttt{11010111}&18&\texttt{01001110}&$r_{41}r_{51}$&$0$\\
21&\texttt{00100100}&16&\texttt{01101101}&$r_{19}$&$0$\\
22&\texttt{01000111}&11&\texttt{01110110}&$r_{24}$&$0$\\
23&\texttt{10111011}&7&\texttt{01110101}&$r_{31}$&$0$\\
24&\texttt{10100100}&12&\texttt{10010101}&$r_{24}$&$0$\\\bottomrule
\end{tabular}
\caption{Leading terms for the cap. The vertices $x_j$ have $L_v\simeq S^3$, and the vertices $y_i$ have $L_v\simeq S^2$. The word $w_j$ is read from $x_j$: the first letter is the first facet crossed. The last column lists the integers $u_i$ of Proposition~\ref{prop:capcert}.}\label{tab:pivots}
\end{table}

\emph{The core.} The vertices $x_j$ and $y_i$ and the basis elements $\mathbf e_j$ and $\mathbf f_i$ for the core, $1\le i,j\le214$, are listed in the ancillary files, where the indices start from $0$. Table~\ref{tab:basis} lists the basis elements $\mathbf e'_j$ of Section~\ref{sec:core} that differ from $\mathbf e_j$, and Table~\ref{tab:corepivots} lists the data of Proposition~\ref{prop:corecert}.

\begin{table}[ht]
\centering\small
\begin{tabular}{rcl}
\toprule
$j$&$x_j$&$\mathbf e'_j-\mathbf e_j$\\\midrule
6&\texttt{01000000}&$\mathbf e_{4}+\mathbf e_{5}$\\
25&\texttt{01101000}&$-\mathbf e_{24}$\\
26&\texttt{01101000}&$-\mathbf e_{25}$\\
31&\texttt{00111000}&$-\mathbf e_{30}$\\
58&\texttt{01011100}&$\mathbf e_{57}$\\
87&\texttt{01111010}&$\mathbf e_{85}$\\
92&\texttt{01000110}&$\mathbf e_{90}$\\
95&\texttt{00100110}&$-\mathbf e_{93}$\\
108&\texttt{01001110}&$\gamma[23]\,\mathbf e_{210}$\\
111&\texttt{00101110}&$\gamma[80,83]\,\mathbf e_{81}$\\
125&\texttt{01010001}&$\mathbf e_{123}+\mathbf e_{124}$\\
130&\texttt{00001001}&$-\gamma[83,110]\,\mathbf e_{127}$\\
139&\texttt{01000101}&$-\gamma[22,50]\,\mathbf e_{90}$\\
146&\texttt{01010101}&$-\mathbf e_{144}$\\
149&\texttt{00001101}&$-\gamma[44,114]\,\mathbf e_{25}+\gamma[44,114]\,\mathbf e_{26}$\\
150&\texttt{01001101}&$\gamma[104,110]\,\mathbf e_{66}$\\
152&\texttt{01001101}&$\gamma[83,109]\,\mathbf e_{22}$\\
153&\texttt{01001101}&$\gamma[83,109]\,\mathbf e_{22}+\mathbf e_{152}$\\
157&\texttt{00011101}&$-\gamma[110,111]\,\mathbf e_{66}$\\
163&\texttt{00000011}&$-\gamma[46,83]\,\mathbf e_{131}+\gamma[46,27,83,25]\,\mathbf e_{164}$\\
164&\texttt{00000011}&$-\gamma[25,27]\,\mathbf e_{130}-\gamma[25,27]\,\mathbf e_{131}$\\
166&\texttt{01000011}&$-\gamma[66,97]\,\mathbf e_{9}$\\
168&\texttt{00100011}&$\gamma[52,67]\,\mathbf e_{33}$\\
170&\texttt{01100011}&$-\gamma[66,109]\,\mathbf e_{144}+\gamma[66,109]\,\mathbf e_{146}+\gamma[66,114]\,\mathbf e_{147}$\\
174&\texttt{01010011}&$\gamma[92,27]\,\mathbf e_{164}-\mathbf e_{173}$\\
181&\texttt{00001011}&$-\mathbf e_{178}$\\
193&\texttt{01111011}&$\gamma[20,23]\,\mathbf e_{176}$\\
194&\texttt{01111011}&$-\gamma[95,77,79,83]\,\mathbf e_{23}+\gamma[95,77]\,\mathbf e_{81}+\gamma[95,77,80,83]\,\mathbf e_{111}$\\
195&\texttt{01111011}&$-\mathbf e_{193}$\\\bottomrule
\end{tabular}
\caption{The basis elements $\mathbf e'_j$ of $\widehat H_4$ for the core that differ from $\mathbf e_j$. Here $\gamma[F_1,\dots,F_r]$ stands for $\gamma_p$, where $p=(F_1,\dots,F_r)$ is the path from $x_j$ that crosses the facets $E_{F_1},\dots,E_{F_r}$ in this order.}\label{tab:basis}
\end{table}

\begin{table}[p]
\centering\scriptsize\setlength{\tabcolsep}{3.2pt}
\begin{tabular}{rrlr|rrlr|rrlr|rrlr}
\toprule
$j$&$\beta(j)$&$w_j$&$u_{\beta(j)}$&$j$&$\beta(j)$&$w_j$&$u_{\beta(j)}$&$j$&$\beta(j)$&$w_j$&$u_{\beta(j)}$&$j$&$\beta(j)$&$w_j$&$u_{\beta(j)}$\\\midrule
1&40&$48$&$-15$&55&162&$104$&$-7$&109&63&$92$&$-12$&163&189&$49$&$-11$\\
2&149&$9$&$-14$&56&31&$55$&$-13$&110&134&$41$&$-19$&164&26&$23,25,27$&$-14$\\
3&100&$97$&$-9$&57&76&$97$&$-13$&111&202&$77,80,83$&$-15$&165&107&$118$&$-6$\\
4&66&$47$&$-14$&58&191&$7$&$-20$&112&123&$28,110$&$0$&166&34&$120$&$-19$\\
5&88&$44$&$-12$&59&156&$21$&$-2$&113&14&$110$&$0$&167&12&$65$&$-22$\\
6&185&$69$&$-20$&60&121&$77$&$-19$&114&148&$58$&$-3$&168&138&$94$&$-18$\\
7&176&$65$&$-17$&61&93&$72$&$-12$&115&94&$89$&$-6$&169&154&$80,32$&$-11$\\
8&135&$64$&$-14$&62&132&$58$&$-16$&116&158&$13$&$-21$&170&140&$94$&$-6$\\
9&101&$97$&$-17$&63&211&$39$&$-6$&117&68&$81$&$-10$&171&117&$57$&$-20$\\
10&159&$39$&$-12$&64&108&$35$&$-14$&118&171&$70$&$-8$&172&13&$18$&$-12$\\
11&112&$17$&$-15$&65&125&$109$&$-5$&119&59&$39$&$-19$&173&8&$62$&$-10$\\
12&168&$62$&$-10$&66&50&$110$&$-1$&120&197&$44$&$-7$&174&143&$83$&$-21$\\
13&192&$69$&$-8$&67&103&$72$&$-19$&121&60&$50$&$-7$&175&16&$13$&$-17$\\
14&126&$23$&$-11$&68&198&$66$&$-14$&122&6&$106$&$-12$&176&77&$20$&$-9$\\
15&41&$72$&$-6$&69&147&$7$&$-16$&123&144&$5$&$-18$&177&109&$35,106$&$-2$\\
16&175&$114$&$-11$&70&102&$5$&$-20$&124&145&$5$&$-20$&178&212&$48$&$-18$\\
17&204&$68$&$-4$&71&61&$60$&$-15$&125&42&$12,27$&$-10$&179&155&$33$&$-16$\\
18&56&$48$&$-19$&72&213&$9$&$-1$&126&11&$104$&$-12$&180&207&$100$&$-18$\\
19&36&$35$&$-2$&73&44&$110$&$-8$&127&209&$110$&$-11$&181&153&$44$&$-20$\\
20&208&$120$&$-14$&74&80&$3$&$-18$&128&75&$87$&$-8$&182&24&$13$&$-11$\\
21&106&$17$&$-9$&75&51&$33$&$-12$&129&160&$48$&$-11$&183&205&$32$&$-15$\\
22&181&$109$&$-8$&76&127&$40$&$-10$&130&133&$83,49,110$&$-5$&184&201&$35$&$-18$\\
23&85&$80$&$-12$&77&157&$7$&$-14$&131&214&$25$&$0$&185&37&$74$&$-14$\\
24&111&$44$&$-5$&78&161&$59$&$-14$&132&62&$69$&$-5$&186&206&$32$&$-9$\\
25&45&$32$&$-6$&79&195&$104$&$-5$&133&128&$78$&$-9$&187&43&$42$&$-10$\\
26&183&$114$&$-6$&80&55&$83$&$-12$&134&110&$86$&$0$&188&89&$50$&$-11$\\
27&141&$77$&$-19$&81&25&$80$&$-10$&135&71&$69$&$-10$&189&165&$78$&$-4$\\
28&96&$110$&$-7$&82&92&$100$&$-15$&136&49&$11$&$-12$&190&57&$120$&$-3$\\
29&151&$116$&$-6$&83&99&$92$&$-14$&137&210&$30$&$-25$&191&142&$17$&$-7$\\
30&58&$72$&$-1$&84&33&$110$&$-8$&138&90&$109$&$-16$&192&52&$52$&$-13$\\
31&35&$100$&$-15$&85&23&$47$&$-4$&139&122&$22,50,52$&$-17$&193&167&$20,23,60$&$-17$\\
32&82&$17$&$-6$&86&194&$28$&$-11$&140&97&$114$&$-15$&194&196&$27$&$-17$\\
33&173&$52$&$-10$&87&39&$103$&$-7$&141&169&$33$&$-14$&195&54&$68$&$-12$\\
34&166&$7$&$-11$&88&5&$83$&$-7$&142&190&$110$&$-12$&196&120&$83$&$-20$\\
35&30&$27$&$-7$&89&83&$27$&$-4$&143&174&$44$&$-7$&197&129&$29$&$-7$\\
36&19&$92$&$-15$&90&28&$22$&$-4$&144&172&$83$&$-14$&198&70&$61$&$-11$\\
37&46&$3$&$-4$&91&203&$21$&$-4$&145&113&$40$&$-6$&199&69&$61$&$-11$\\
38&186&$53$&$0$&92&188&$77$&$-12$&146&105&$114$&$-10$&200&98&$111$&$-8$\\
39&87&$24$&$-13$&93&118&$49$&$-14$&147&4&$66$&$-6$&201&84&$50$&$-8$\\
40&150&$115$&$-19$&94&27&$17$&$-14$&148&115&$69$&$-18$&202&184&$92$&$-7$\\
41&1&$79$&$-13$&95&114&$38$&$-18$&149&3&$118$&$-4$&203&17&$59$&$-18$\\
42&73&$17$&$-18$&96&199&$23$&$-16$&150&48&$104$&$-10$&204&91&$106$&$-17$\\
43&146&$23$&$-19$&97&139&$13$&$-15$&151&29&$11$&$-15$&205&182&$95$&$-8$\\
44&22&$100$&$-18$&98&200&$16$&$-12$&152&179&$83$&$-8$&206&47&$18$&$-14$\\
45&86&$98$&$-8$&99&2&$30$&$-12$&153&74&$120$&$-6$&207&178&$27$&$-6$\\
46&119&$66$&$-5$&100&163&$11$&$-5$&154&78&$41$&$-10$&208&18&$7$&$-8$\\
47&38&$124$&$-13$&101&7&$30$&$-15$&155&170&$95,47$&$0$&209&130&$44$&$-8$\\
48&65&$17$&$-21$&102&67&$55$&$-10$&156&131&$72$&$-15$&210&137&$97$&$0$\\
49&152&$23$&$-22$&103&79&$92$&$-10$&157&116&$114$&$-5$&211&64&$88$&$-13$\\
50&136&$116$&$-9$&104&32&$36$&$-11$&158&124&$46$&$-3$&212&180&$79$&$-9$\\
51&95&$82$&$-9$&105&187&$57$&$-15$&159&104&$69$&$-18$&213&164&$100$&$-1$\\
52&20&$72$&$-8$&106&15&$44$&$-16$&160&10&$88$&$-8$&214&72&$118$&$-24$\\
53&21&$48$&$-16$&107&81&$37$&$-11$&161&9&$50$&$-14$&&&&\\
54&193&$59$&$-9$&108&177&$21,92$&$-2$&162&53&$23$&$-13$&&&&\\\bottomrule
\end{tabular}
\caption{Leading terms for the core. A path $w_j=(F_1,\dots,F_r)$ from $x_j$ to $y_{\beta(j)}$ is written $F_1,\dots,F_r$. The last column of each block lists the integer $u_i$ of Proposition~\ref{prop:corecert} for $i=\beta(j)$.}\label{tab:corepivots}
\end{table}

\section{Algorithmic verification}\label{app:computations}

All computations are exact over $\Z$, $\Z[\sqrt2]$ and $\F_2$, using the data of Section~\ref{sec:data} and Appendix~\ref{app:data}. The data and verification scripts are included in the ancillary files.

\subsection{Polytopes, color spans, and compatibility}
We compute all inner products $\langle e_i,e_j\rangle$ and $\langle g_i,g_j\rangle$ and enumerate all nonempty cliques of the two nerves ($33{,}644$ for $P^7$ and $207{,}902$ for $Q$). The ideal vertices are found as the null vectors $e_i+e_j$ (resp.\ $g_i+g_j$) with inner product $-1$. This checks Lemmas~\ref{lem:P7} and~\ref{lem:Q}, and we also check that $\jmath$ permutes $\{g_i\}$. By Gaussian elimination over $\F_2$, the colors of $I\setminus\{8\}$ span $V$ and those of $I_{(8)}$ span $U$ for $P^7$, while the colors of $\If$, $I_{(0)}$, and $I_{(181)}$ each span $U$ for $Q$. Every tangent pair of $P^7$ and every tangent fixed pair of $Q$ has distinct colors. Lemma~\ref{lem:compatible} is verified directly on $\{F_j:j\sim8\}$ and $\{E_1,\dots,E_{27}\}$.

\subsection{Proposition~\ref{prop:links} and Lemma~\ref{lem:peripheral}}
For each subcomplex, explicit generator-eliminating Tietze transformations reduce $\pi_1$ to the trivial group ($41{,}472$ steps across the $256$ links $L_v$ for $P^7$, and $202{,}752$ for $Q$). For homology, we cancel pairs $(\sigma,\tau)$ as in Section~\ref{sec:local} for $k=1,\dots,7$ in turn. At each step we pick a pair that minimizes $(r-1)(c-1)$, where $r$ and $c$ are the row and column supports of $\tau$ and $\sigma$ in the current differential, with a fixed tie-breaking rule. Afterwards, the identities \eqref{eq:sdr} are checked directly against the original boundary matrices on all basis elements $e_{v,S}$ ($468{,}224$ for $P^7$ and $2{,}710{,}016$ for $Q$).

For Lemma~\ref{lem:peripheral}(1)--(2), whether $|\{\alpha,\beta\}\cap D(v)|=1$ depends only on $v\cdot(c_\alpha+c_\beta)$. At each ideal vertex, the vector $\delta=c_\alpha+c_\beta$ is the same for all relevant tangent pairs, so it suffices to check both values of $v\cdot\delta$ at every ideal vertex. Part~(3) is verified from the fundamental cycle periods of a spanning tree of the Schreier graph on $V$.

\subsection{Propositions~\ref{prop:capcert} and~\ref{prop:corecert}}\label{app:core}
For the cap $C$, we evaluate $D^{(1)}$ ($16$ terms) and $D^{(2)}$ ($163$ terms) exactly in $\Z\Gamma_C$, comparing descending words modulo commutation of adjacent letters. Apart from the $24$ leading terms there are $155$ terms, and the minimum gap $a(w)-a(w_j)-u_i+u_{\beta(j)}$ is $1$.

For the core $X$, we evaluate $D^{(1)}$ ($1166$ terms) and $D^{(2)}$ ($10947$ terms) in the same way. All coefficients are $\pm1$, except for $80$ terms with coefficients $\pm2$. The truncation $D'^{\,\le2}$ is assembled via \eqref{eq:Dprime} and Table~\ref{tab:basis}. Because paths $pw$ need not be descending, group elements are reduced to normal form in the right-angled Coxeter group $W$ by cancelling equal letters separated only by commuting letters and taking the lexicographically smallest commutation-equivalent word. Collecting equal group elements yields $13074$ nonzero terms in $D'^{\,\le2}$. Among the $214$ paths $w_j$, $203$ have length $1$, $6$ have length $2$, and $5$ have length $3$. The $12860$ non-leading terms have minimum gap $a(q)-a(w_j)-u_i+u_{\beta(j)}=1$, the maximum of $a(w_j)-a(p)+u_i-u_{\beta(j)}$ in Proposition~\ref{prop:corecert}(3) is $50$, and the longest chains in $\prec$ have length $3$.

Both Propositions~\ref{prop:capcert} and~\ref{prop:corecert} were verified by two independent programs, the second comparing group elements via exact matrix multiplication of the reflections $r_i$ over $\Z[\sqrt2]$. The operators $D^{(1)}$ and $D^{(2)}$ were likewise computed by two independent programs. The weights, basis $(\mathbf e'_j)$, and leading-term data were found via linear and integer programming, and the proofs rely only on Lemmas~\ref{lem:compatible} and~\ref{lem:peripheral} and Propositions~\ref{prop:capcert} and~\ref{prop:corecert}.


\clearpage
\begin{thebibliography}{IMM24}

\bibitem[Ago13]{Agol13}
I.~Agol, The virtual Haken conjecture, Doc. Math. 18 (2013), 1045--1087. With an appendix by I.~Agol, D.~Groves, and J.~Manning.

\bibitem[ALR01]{ALR01}
I.~Agol, D.~D. Long, and A.~W. Reid, The Bianchi groups are separable on geometrically finite subgroups, Ann. of Math. (2) 153 (2001), no.~3, 599--621.

\bibitem[BL12]{BL12}
A.~Bartels and W.~L\"uck, The Borel conjecture for hyperbolic and CAT(0)-groups, Ann. of Math. (2) 175 (2012), no.~2, 631--689.

\bibitem[BB97]{BB97}
M.~Bestvina and N.~Brady, Morse theory and finiteness properties of groups, Invent. Math. 129 (1997), no.~3, 445--470.

\bibitem[Bro65]{Brown65}
R.~Brown, The twisted Eilenberg--Zilber theorem, Simposio di Topologia (Messina, 1964), Edizioni Oderisi, Gubbio, 1965, pp.~33--37.

\bibitem[Cra04]{Crainic04}
M.~Crainic, On the perturbation lemma, and deformations, preprint (2004), arXiv:math/0403266.

\bibitem[Dav08]{Davis08}
M.~W. Davis, The geometry and topology of Coxeter groups, London Mathematical Society Monographs Series, vol.~32, Princeton University Press, Princeton, NJ, 2008.

\bibitem[Duf10]{Dufour10}
G.~Dufour, Notes on right-angled Coxeter polyhedra in hyperbolic spaces, Geom. Dedicata 147 (2010), 277--282.

\bibitem[ERT12]{ERT12}
B.~Everitt, J.~G. Ratcliffe, and S.~T. Tschantz, Right-angled Coxeter polytopes, hyperbolic six-manifolds, and a problem of Siegel, Math. Ann. 354 (2012), no.~3, 871--905.

\bibitem[Far71]{Farrell71}
F.~T. Farrell, The obstruction to fibering a manifold over a circle, Actes du Congr\`es International des Math\'ematiciens (Nice, 1970), Tome~2, Gauthier-Villars, Paris, 1971, pp.~69--72.

\bibitem[Far72]{Farrell72}
F.~T. Farrell, The obstruction to fibering a manifold over a circle, Indiana Univ. Math. J. 21 (1971/72), 315--346.

\bibitem[Fis23]{Fisher23}
S.~P. Fisher, Algebraic fibring of a hyperbolic $7$-manifold, Bull. Lond. Math. Soc. 55 (2023), no.~3, 1347--1357.

\bibitem[GL14]{GL14}
T.~Gelander and A.~Levit, Counting commensurability classes of hyperbolic manifolds, Geom. Funct. Anal. 24 (2014), no.~5, 1431--1447.

\bibitem[GPS88]{GPS88}
M.~Gromov and I.~Piatetski-Shapiro, Non-arithmetic groups in Lobachevsky spaces, Inst. Hautes \'Etudes Sci. Publ. Math. 66 (1988), 93--103.

\bibitem[IMM23]{IMM23}
G.~Italiano, B.~Martelli, and M.~Migliorini, Hyperbolic $5$-manifolds that fiber over $S^1$, Invent. Math. 231 (2023), no.~1, 1--38.

\bibitem[IMM24]{IMM24}
G.~Italiano, B.~Martelli, and M.~Migliorini, Hyperbolic manifolds that fibre algebraically up to dimension $8$, J. Inst. Math. Jussieu 23 (2024), no.~2, 609--646.

\bibitem[IMN26]{IMN26}
G.~Italiano, M.~Migliorini, and A.~Ng, Improved algebraic fibrations of high-dimensional hyperbolic groups, preprint (2026), arXiv:2606.05091.

\bibitem[JNW21]{JNW21}
K.~Jankiewicz, S.~Norin, and D.~T. Wise, Virtually fibering right-angled Coxeter groups, J. Inst. Math. Jussieu 20 (2021), no.~3, 957--987.

\bibitem[Kie25]{Kielak25}
D.~Kielak, Virtual fibring of manifolds and groups, preprint (2025), arXiv:2510.01805.

\bibitem[Kol26]{Kolpakov26}
A.~Kolpakov, Sharp dimension bounds for right-angled hyperbolic polytopes, preprint (2026), arXiv:2610.38900.


\bibitem[PV05]{PV05}
L.~Potyagailo and E.~Vinberg, On right-angled reflection groups in hyperbolic spaces, Comment. Math. Helv. 80 (2005), no.~1, 63--73.

\bibitem[Ran95]{Ranicki95}
A.~Ranicki, Finite domination and Novikov rings, Topology 34 (1995), no.~3, 619--632.

\bibitem[Sch06]{Schuetz06}
D.~Sch\"utz, Finite domination, Novikov homology and nonsingular closed $1$-forms, Math. Z. 252 (2006), no.~3, 623--654.

\bibitem[Sk\"o06]{Skoldberg06}
E.~Sk\"oldberg, Morse theory from an algebraic viewpoint, Trans. Amer. Math. Soc. 358 (2006), no.~1, 115--129.

\bibitem[Vin85]{Vinberg85}
\`E.~B. Vinberg, Hyperbolic reflection groups, Russian Math. Surveys 40 (1985), no.~1, 31--75.

\bibitem[Wal99]{Wall99}
C.~T.~C. Wall, Surgery on compact manifolds, 2nd ed., Mathematical Surveys and Monographs, vol.~69, American Mathematical Society, Providence, RI, 1999.

\bibitem[Wis21]{Wise21}
D.~T. Wise, The structure of groups with a quasiconvex hierarchy, Annals of Mathematics Studies, vol.~209, Princeton University Press, Princeton, NJ, 2021.

\end{thebibliography}
\end{document}